

\font\scX=cmcsc10
\font\rmVIII=cmr8
\font\slVIII=cmsl8

\font\bfXII=cmbx10 scaled \magstep1
\def\statement#1#2{\medbreak\noindent{\bf#1.\enspace}{\sl#2\par}
\ifdim\lastskip<\medskipamount\removelastskip\penalty55 \medskip\fi}
\def\proof{\smallskip \noindent {\bf Proof. }}
\def\square{\hbox{$\sqcap\!\!\!\!\sqcup$}}
\def\endproof{\hfill \square}
\def\cite#1{{\rm[\bf #1\rm]}}
\def\stdbib#1#2#3#4#5#6{\smallskip
\item{[#1]} #2, ``#3,'' {\sl #4} {\bf #5} #6.}
\def\bib#1#2#3#4{\smallskip \item{[#1]} #2, ``#3,'' {#4.}}
\def\stackrel#1#2{\buildrel #1 \over #2}
\magnification=\magstep1
\nopagenumbers
\voffset=2\baselineskip
\advance\vsize by -\voffset
\headline{\ifnum\pageno=1\hfil\else\runninghead\fi}
\def\runninghead{\scX\hfil circle actions and partial
automorphisms\hfil\folio}
\def\cstar{$C^*$}
\def\E{{\cal E}}
\def\H{{\cal H}}
\def\K{{\cal K}}
\def\Z{{\bf Z}}
\def\N{{\bf N}}
\def\S1{S^1}
\def\for{\hbox{\ \ \ for\ }}
\def\th{\null^{\rmVIII th}}
\def\Dom{\hbox{Dom}}
\def\Im{\hbox{Im}}
\def\KK{\hbox{$K\!K$}}
\def\SameAuthor{\underbar{\hbox{\ \ \ \ \ \ \ \ \ \ }}}
\def\MRClass{\rmVIII 1991 \slVIII MR Subject Classification:
\rmVIII Primary 46L05, 46L40, 46L45, 46L55, 46L80, 19K35; Secondary 47B35,
34C35.}
\def\support{\rmVIII Partially supported by FAPESP, Brazil. On leave
from the University of S\~ao Paulo.}
\def\AP{1}
\def\Brown{2}
\def\BGR{3}
\def\CuntzOnI{4}
\def\CuntzOnII{5}
\def\CuntzKK{6}
\def\CuntzNewLook{7}
\def\FellI{8}
\def\FellII{9}
\def\Jensen{10}
\def\KishimotoTakai{11}
\def\Kumjian{12}
\def\PaschkeEndom{13}
\def\PaschkeCircle{14}
\def\Pedersen{15}
\def\Phillips{16}
\def\PV{17}
\def\RieffelInduced{18}
\def\RieffelMorita{19}
\def\RieffelProper{20}

\null
\baselineskip=18pt
\centerline{\bfXII CIRCLE ACTIONS ON C\raise.5ex\hbox{*}-ALGEBRAS, PARTIAL}
\centerline{\bfXII AUTOMORPHISMS AND A GENERALIZED}
\centerline{\bfXII PIMSNER--VOICULESCU EXACT}
\centerline{\bfXII SEQUENCE}

\bigskip
\bigskip
\centerline{\scX Ruy Exel\footnote{$\null^*$}{\support}
\footnote{\null}{\MRClass}}

\bigskip
\baselineskip=12pt
\centerline{\it Department of Mathematics and Statistics}
\centerline{\it University of New Mexico}
\centerline{\it Albuquerque, New Mexico 87131}
\centerline{e-mail: exel@math.unm.edu}

\vfill
\midinsert\narrower\noindent
{\bf Abstract.} We introduce a method to study
C$\null^*$-algebras possessing an action of the circle group, from the
point of view of its internal structure and its K-theory.  Under
relatively mild conditions our structure Theorem shows that any
C$\null^*$-algebra, where an action of the circle is given, arises as
the result of a construction that generalizes crossed products by the
group of integers.

Such a generalized crossed product construction is carried out for
any partial automorphism of a C$\null^*$-algebra, where by a partial
automorphism we mean an isomorphism between two ideals of the given
algebra.

Our second main result is an extension to crossed products by
partial automorphisms, of the celebrated Pimsner-Voiculescu exact
sequence for K-groups.

The representation theory of the algebra arising from our
construction is shown to parallel the representation theory for
C$\null^*$-dynamical systems. In particular, we generalize several of
the main results relating to regular and covariant representations of
crossed products.
\endinsert\vfill\eject

\beginsection 1. Introduction

Given a group action on a manifold $M$, one of the main goals of the
Dynamical Systems specialist is to describe $M$ in terms of the
elementary components singled out by the presence of the action as,
for example, orbits and fixed points.  The simplest case in which such
a description can be thoroughly carried out is that of a free action
of a compact group: $M$ can be then described as a principal bundle
over the quotient space.

While the interaction between Operator Algebras and Dynamical
Systems has been very intense in the last several decades, little has
been accomplished in addressing the above question from the Operator
Algebras point of view.  The \cstar-Dynamical Systems counterpart of
that program would be to describe the structure of a \cstar-algebra,
where a group action is given, in terms of elementary data which,
supposedly should be extracted from the action.

Among the few cases in which the above mentioned task was
successfully carried out, is a result by Paschke (Theorem 2.3 in
\cite{\PaschkeCircle}) in which it is shown that, under certain
circumstances, a \cstar-algebra carrying an action of the circle group
can be described as the crossed product of its fixed point subalgebra
(the counterpart of the quotient space) by an action of the integers.
The main hypothesis in that Theorem is that the action have ``large
spectral subspaces'' (see \cite{\PaschkeCircle} and
\cite{\KishimotoTakai}), a condition that, at least in the case of the
circle, singles out the non-commutative analog of free actions.

If one goes as far as to accept that a crossed product is the
non-commutative version of a principal bundle, then Paschke's result
can be considered as an ipsis-literis generalization of the well known
fact about free actions of compact groups mentioned above.

The circle of ideas around the notion of ``large spectral
subspaces'' has now a long history.  Without attempting a
comprehensive account, we should mention that a version of this notion
appeared in a paper by Fell \cite{\FellI} in 1969 under the name of
homogeneous actions, later renamed by Fell himself \cite{\FellII} to
saturated actions.  Related notions were also studied by Phillips
\cite{\Phillips}. Recent work of Rieffel \cite{\RieffelProper}
contains yet another version of this concept, providing a far reaching
generalization of the Takesaki-Takai duality (Corollary 1.7 of
\cite{\RieffelProper}, see also \cite{\KishimotoTakai}).

But, in the same way that free actions are not the rule, the
conditions so far alluded to, exclude many actions of $\S1$ which one
would still like to investigate.

The purpose of the present work is to introduce a method which
allows for a description of the structure of \cstar-algebras carrying
circle actions which are not supposed to have large spectral
subspaces.  While our method does not include all possible circle
actions, since we assume our actions to be semi-saturated (see below),
the gain in generality is very significant in the sense that a wealth
of new examples becomes tractable using our theory.  A typical such
example is the action of $\S1$ on the Toeplitz algebra, given by
conjugation by the diagonal unitaries diag$(1, z, z^2, ...)$, for $z
\in \S1$.  The condition of having large spectral subspaces fails for
this action.

Our method consists in first, introducing a construction, inspired
on the crossed product construction, which produces a \cstar-algebra
equipped with a circle action.  We then go on to show that any
\cstar-algebra having an action of $\S1$ arises as the result of our
construction as long as the action satisfies some relatively mild
restrictions. In other words we provide a means to ``disassemble'' a
\cstar-algebra possessing a circle action, reverting, in a sense, the
crossed product construction.

While a crossed product by $\Z$ depends on an automorphism of the
given \cstar-algebra $A$ (i.e., an action of $\Z$), our construction
requires a partial automorphism. Precisely speaking, a partial
automorphism of $A$ is a triple $\Theta = (\theta, I, J)$ where $I$
and $J$ are ideals in $A$ and \ $\theta \colon  I \rightarrow J$ \ is a
\cstar-algebra isomorphism. Given a partial automorphism, we construct
its ``covariance algebra'' which we denote by $C^*(A,\Theta)$.  When
both ideals agree with $A$, then our construction becomes the usual
crossed product construction, that is, $C^*(A,\Theta)$ becomes
$A\times_\theta \Z$. For that reason, $C^*(A,\Theta)$ should be
considered as the crossed product of $A$ by the partial automorphism
$\Theta$.

Our second major objective, accomplished by Theorem (7.1), is a
generalization of the celebrated Pimsner--Voiculescu exact sequence
\cite{\PV}.  Precisely, we get the following exact sequence of
$K$-groups
$$ \matrix {
K_0(J) &\buildrel {i_* - \theta^{-1}_*} \over \longrightarrow
&K_0(A) &\buildrel {i_*} \over \longrightarrow
&K_0(C^*(A,\Theta)) \cr
\cr
\uparrow & & & &\downarrow \cr
\cr
K_1(C^*(A,\Theta)) &\buildrel {i_*} \over \longleftarrow
&K_1(A) &\buildrel {i_* - \theta^{-1}_*} \over \longleftarrow
&K_1(J)
}$$

As in most proofs of the Pimsner--Voiculescu exact sequence and
related results (\cite{\CuntzOnI}, \cite{\CuntzOnII}, \cite{\CuntzKK},
\cite {\AP}), we derive our exact sequence from the $K$-theory exact
sequence for a suitable Toeplitz extension. The crucial step, as it is
often the case, is to show that $A$ has the same $K$-groups as the
Toeplitz algebra. We do so by showing that these algebras are, in
fact, \KK-equivalent.

My first attempt at proving the exactness of the sequence above was,
of course, by trying to deduce it from the well known result of
Pimsner and Voiculescu. After failing in doing so, I am now tempted to
believe that this cannot be done. Our proof is done from scratch and,
given a rather involving use of \KK-theory, it turns out considerably
longer then the available proofs of the original result.

We would like to thank Bill Paschke for bringing to our attention
his paper \cite{\AP} with J. Anderson, where they generalize a
result, from unpublished lecture material of Arveson's, as well as
from Proposition (5.5) in \cite{\CuntzKK}, from which the crucial step
in \cite{\PV} follows. Our generalization of these ideas plays a
central role in the proof of our result.

One intersting aspect, central in our use of \KK-theory, is worth
mentioning here. If $A$ and $B$ are \cstar-algebras, then $\KK(A,B)$
may be described, as was shown by Cuntz \cite{\CuntzNewLook}, by the
set of homotopy classes of homomorphisms from $qA$ to the multiplier
algebra of $B\otimes K$. Nevertheless the \KK-theory elements that we
need to introduce have no easy description in such terms.  Instead we
exhibit these elements by replacing $B\otimes K$, above, by an algebra
which contains $B$ as a full corner and hence is stably isomorphic to
$B$ (at least in the separable case).

We feel that our structure Theorem, used in conjunction with the
generalized Pimsner--Voiculescu exact sequence above, adds to the
growing collection of powerful tools, developed since the early
eighties, designed to compute $K$-groups for \cstar-algebras.

Several attempts to generalize crossed products have been made by
many authors. Among those we should mention Kumjian's work on
Localizations \cite{\Kumjian}.  Another well known example is the
theory of crossed products by endomorphisms developed by Paschke in
\cite{\PaschkeEndom}. Paschke's crossed products carry a circle
action, which can easily be made to fall under our restrictions. So,
in a certain sense, his algebras are special cases of our
construction.  In addition, we have recently learned that L.~Brown has
developed a theory which formalizes the concept of crossed products by
imprimitivity bimodules.

Our construction should be regarded as a generalization of crossed
products, only as long as the group $\Z$, of integers, is concerned.
With some more work, it appears to me that one could attempt to widen
the present methods to include a larger class of groups.  It would be
largely desirable, although quite likely very difficult, to be able to
study along the present lines, actions of non-compact groups such as
the group of real numbers.
Florin Pop pointed out to me recently that the definition of a
``partial action'' of a discrete group, as well as that of the
corresponding covariance algebra, could be obtained by trivial
modifications of our definitions.

As we already indicated, our motivation, rather than to produce new
classes of \cstar-algebras, is to attempt a description of the
structure of \cstar-algebras possessing a circular symmetry,
represented by an action of the circle.  A partial automorphism arises
from a given action $\alpha$ of $\S1$ on a \cstar-algebra $B$, in the
following way.  First, one lets $B_1$ be the first spectral subspace
of $\alpha$, i.e., $B_1 = \{b \in B: \alpha_z (b) = zb, z\in \S1\}$.
The set $B^*_1 B_1$ (meaning, according to convention (2.2) adopted
throughout this paper, the closed linear span of the set of products)
is an ideal of the subalgebra $B_0$ of fixed points for $\alpha$ and
the same is true with respect to $B_1 B_1^*$.  These ideals are
strongly Morita equivalent \cite{\RieffelInduced},
\cite{\RieffelMorita} with $B_1$ playing the role of the imprimitivity
bimodule.  Morita equivalent \cstar-algebras are quite often
isomorphic to each other and indeed, under the assumption that the
algebras be stable with strictly positive elements, they are forcibly
isomorphic \cite{\BGR}.  The two isomorphic ideals are thus the
ingredients of our partial automorphism.

Our structure Theorem (4.21), states that $B$ is (isomorphic to) the
covariance algebra of that partial automorphism as long as $B$ is
generated, as a \cstar-algebra, by the union of $B_0$ and $B_1$.  This
condition, which we call semi-saturation, is a weakening of the
condition of having large spectral subspace and is the point of
departure for our theory.  As opposed to what happens to the latter,
even when absent, our condition can be forced upon the action by
restricting one's attention to the subalgebra of $B$ that $B_0 \cup
B_1$ happen to generate.

The algebraic formalism developed here seems to flow with such a
naturality that it is perhaps a bit surprising that it has been
overlooked until now.  Nevertheless, it does not seem possible to
extend our methods to rings not possessing a \cstar-algebra structure
since we make extensive use of the existence of approximate units and
facts like ideals of ideals of a \cstar-algebra are, themselves,
ideals of that algebra, or that the intersection of two ideals equals
their product.

Permeating most of our techniques is a concept we have not tried to
formalize, but I think the effort to do so seems worthwhile.  The
reader is invited to compare the definition of multipliers of
\cstar-algebras on one hand, and (4.11) and (4.13) on the other, and
he will likely see the rudiments of a concept which deserves the name
of partial multipliers.

After a short section (2), intended mainly to fix some notation, we
describe in section (3) our construction of the covariance algebra
associated to a partial automorphism.  The following section, numbered
(4) is where our structure Theorem (4.21), is proved.  The fifth
section deals with the representation theory for our covariance
algebras and to some extent can be regarded as the taming of the
algebraic properties of partial isometric operators on Hilbert's
space.  With surprising ease, partial isometries are made to play the
role usually played by implementing unitaries for representations of
crossed product algebras.

In section (6) we introduce the Toeplitz algebra associated to a
partial automorphism and describe it, within our theory, as a
covariance algebra, as well.  In doing so, we are able to obtain a
crucial universal property, Lemma (6.7), characterizing
representations of the Toeplitz algebra. The seventh and last section
is where our main $K$-theoretical work is developed and where we prove
the existence of the generalized Pimsner--Voiculescu exact sequence
(Theorem 7.1).

Our notation is reasonably standard except, possibly, for the
definitions in (2.2) and (5.4) as well as for our use of the symbol
$\bigoplus$ after (6.3). The reader is advised to go over these
immediately, in order to avoid possible surprises.

This work is an extended version of a paper by the author, entitled
``The Structure of Actions of the Circle Group on \cstar-Algebras''.

\beginsection 2. Spectral Subspaces

This section is concerned with some preliminaries about
\cstar-dynamical systems based on the circle group.  Let $B$ be a
fixed \cstar-algebra and $\alpha$ an action of $\S1$ on $B$.

\statement{2.1. Definition}{For each $n\in \Z$ the $n\th$ spectral
subspace for $\alpha$ is defined by
$$ B_n = \{b\in B: \alpha_z (b) = z^nb \for z\in \S1\}.$$}

\smallskip It is an easy matter to verify that $B_n B_m \subseteq
B_{n+m}$ and that $B_n^* = B_{-n}$. Regarding the product $B_n B_m$
just mentioned, we adopt the following convention:

\statement{2.2. Definition}{If $X$ and $Y$ are subsets of a
\cstar-algebra then $XY$ denotes the \underbar{closed} linear span of
the set of products $xy$ with $x \in X$ and $y \in Y$.}

\smallskip A simple fact we shall make extensive use of, is the
following.

\statement{2.3. Proposition}{For each $n \in \Z$ one has that $B^*_n
B_n$ is a closed two sided ideal of the fixed point subalgebra $B_0$.}

\statement{2.4. Definition}{The $n\th$ spectral projection for the
action $\alpha$ is the transformation
$$ P_n \colon  B \rightarrow B$$
defined by
$$ P_n (b) = \int_{S^1} z^{-n} \alpha_z (b) dz\ ,\ b\in B.$$}

\smallskip It is well known that $P_n$ is a contractive projection whose
image is $B_n$.

\statement{2.5. Proposition} Let $b \in B$ and $\phi$ be a
continuous linear functional on $B$.
\medskip
\itemitem{(a)} If $\phi(P_n(b)) = 0$ for all $n \in \Z$, then $\phi(b)=0$.
\medskip
\itemitem{(b)} If $P_n(b) = 0$ for all $n \in \Z$, then $b=0$.
\medskip
\itemitem{(c)} $\bigoplus_{n\in \Z} B_n$ is dense in $B$.

\proof (cf. \cite{\PaschkeCircle}, Proposition 2.1).  For the first
statement it is enough to note that the $n\th$ Fourier coefficient of
$z \in \S1 \rightarrow \phi (\alpha_z (b))$ is given by
$\phi(P_n(b))$. Hence, if $\phi(P_n(b)) = 0$ for all $n$ we have that
$\phi (\alpha_z (b))$ is identically zero, as a function of $z$ and,
in particular, $\phi(b) = 0$. Finally (b) and (c) follow from (a) and
the Hahn-Banach Theorem. \endproof

\statement{2.6. Proposition} {If $(e_\lambda )_{\lambda \in \Lambda}
$ is an approximate identity for $B_n^* B_n$ then, for each $x \in
B_n$, we have $x =\lim_\lambda x e_\lambda$.}

\proof We have
$$ \| x - x e_\lambda \|^2 = \| (x -x e_\lambda )^* (x-x e_\lambda )
\| = \| x^* x - x^* x e_\lambda - e_\lambda x^* x + e_\lambda x^* x
e_\lambda \| \leq$$
$$ \| x^* x - x^* x e_\lambda \| + \| e_\lambda \| \| x^* x - x^* x
e_\lambda \| \rightarrow 0.$$
\endproof

\smallskip An immediate consequence is:

\statement{2.7. Corollary}{For each $n \in \Z$ one has $B_n B_n^*
B_n = B_n$.}

\smallskip Since the product of two ideals in a \cstar-algebra equals
their intersection we have:

\statement{2.8. Proposition} {If $n$ and $m$ are integers then
$B_n^* B_n B^*_m B_m = B_m^* B_m B_n^* B_n$.}

\smallskip The following concludes our preparations.

\statement{2.9. Proposition}{Let $B$ and $B^\prime$ be
\cstar-algebras and let $\alpha$ and $\alpha^\prime$ be actions of
$\S1$ on $B$ and $B^\prime$, respectively. Suppose $\psi \colon  B
\rightarrow B^\prime$ is a covariant homomorphism. If the restriction
of $\psi$ to the fixed point subalgebra $B_0$ is injective, then
$\psi$ itself, is injective.}

\proof Assume $b \in B$ is such that $\psi(b) = 0$. Then, denoting
by $P_n$ and $P_n^\prime$ the respective spectral projections, we have
$$\psi(P_n(b)P_n(b)^*) = \psi(P_n(b))\psi(P_n(b))^* =
P_n^\prime(\psi(b)) P_n^\prime(\psi(b))^* = 0.$$

If one now notes that $P_n(b)P_n(b)^* \in B_0$, the hypothesis is
seen to imply that $P_n(b) = 0$ and so, by (2.5.b), we have $b=0$.
\endproof

\beginsection 3. Partial Automorphisms and their Covariance Algebras

In this section we describe a generalization of the concept of crossed
products by an automorphism.  For that purpose, let $A$ be a
\cstar-algebra considered fixed throughout the present section.

\statement{3.1. Definition}{A partial automorphism of $A$ is a
triple $\Theta = (\theta, I, J)$ where $I$ and $J$ are ideals in $A$
(always assumed closed and two sided) and \ $\theta \colon  I \rightarrow J$
\ is a \cstar-algebra isomorphism.}

\smallskip If such a partial automorphism is given, we let, for each
integer $n, D_n$ denote the domain of $\theta^{-n}$ with the
convention that $D_0 = A$ and $\theta^0 $ is the identity automorphism
of $A$.  The domain of $\theta^{-n} $ is clearly the image of
$\theta^n$ so this provides an equivalent definition of $D_n$.

Alternatively, we can give an inductive definition for these objects
by letting $D_0 = A$,
$$ D_{n+1} = \{a\in J: \theta^{-1} (a) \in D_n\}$$
for $n \geq 0$ and
$$ D_{n-1} = \{a \in I: \theta (a) \in D_n\}$$
for $n \leq 0$.

According to this, one clearly has $D_1 = J$ and $D_{-1} = I$.  Of
course, unless $I$ and $J$ have a substantial intersection, the sets
$D_n$ would be rather small.  The extreme case in which $I \cap J =
\{0 \}$ will see $D_n$ being the singleton $\{0 \}$ for all $n$,
except for $n = -1, 0, 1$.

One of the simplest examples of partial automorphisms is obtained
when one lets $A = {\bf C}^m$, $I = \{(x_i) \in {\bf C}^m : x_m = 0\},
J = \{(x_i) \in {\bf C}^m: x_1 = 0 \}$ and $\theta$ be the forward
shift
$$ \theta (x_1, \ldots, x_{m-1}, 0) = (0, x_1, \ldots, x_{m-1}).$$

In this case $D_n$ becomes the set of all $m$-tuples having $n$
leading zeros when $n \geq 0$ or $|n|$ trailing zeros if $n \leq 0$.

As it will turn out, the covariance algebra for this example is
isomorphic to the algebra of $m\times m$ complex matrices.

\statement{3.2. Proposition}{For each integer $n$, $D_n$ is an ideal
in $A$.}

\proof By definition the assertion is obvious for $n = -1, 0, 1$.
Arguing by induction assume that $n \geq 0$ and that $D_n$ is and
ideal in $A$.  Then $D_{n+1}$ is clearly an ideal in $J$, being the
inverse image of $D_n$ under $\theta^{-1}$.  Since an ideal of an
ideal of a \cstar-algebra is always an ideal of that \cstar-algebra
(by existence of approximate identities) we have that $D_{n+1}$ is an
ideal in $A$.  A symmetric argument yields the result for negative
values of $n$.
\endproof

\smallskip Let's agree to call by the name of chain any finite sequence
$(a_0, a_1, \ldots, a_n)$ of elements in $A$ such that $a_0 \in I$,
$a_n \in J$ and $a_i \in I \cap J$ for $i = 1, \ldots, n-1$,
satisfying $\theta (a_{i-1} ) = a_i$ for $i = 1, 2, \ldots, n$.  The
integer $n$ will be called the length of said chain.

The concept of chain can be used for giving yet another definition
of $D_n$.  Namely, $D_n$ is the set of elements a in $A$ which admit a
chain of length $|n|$, ending in $a$ in case $n \geq 0$ or, beginning
in $a$ if $n \leq 0$.

The following proposition can be easily proven if one thinks in
terms of chains.  It is nevertheless crucial for what follows.

\statement{3.3. Proposition}{If $n$ and $m$ are integers then
$\theta^{-n} (D_n \cap D_m) \subseteq D_{m-n}$.  In addition, if $x$
is in $D_n \cap D_m$ then \ $\theta^{n-m} (\theta^{-n} (x)) =
\theta^{-m} (x)$.}

\smallskip Denote by $L$ the subspace of $\ell_1 (\Z, A)$ formed by all
summable sequences $(a (n))_{n\in \Z}$ such that $a(n) \in D_n$ for
each $n$.  We propose to equip $L$ with an involutive Banach algebra
structure.  For that purpose we define, for $a$ and $b$ in $L$,
$$ (a * b) (n) = \sum^\infty_{k=-\infty} \theta^k ( \theta^{-k} (a
(k)) b (n-k))$$
$$ (a^*) (n) = \theta^n (a (-n)^*)$$
$$ \| a \| = \sum^\infty_{n=-\infty} \| a (n) \|.$$

We next verify some of the axioms of involutive Banach algebras for
the multiplication, involution and norm defined above.  But, before
that, we should note that our multiplication is well defined since,
for each $k, \theta^{-k} (a(k))$ is in $D_{-k}$ while $b(n-k) $ is in
$D_{n-k}$.  The product $\theta^{-k}(a(k)) b(n-k)$ is therefore in the
intersection $D_{-k} \cap D_{n-k}$.  By Proposition (3.3) it follows
that the $k\th$ summand in our definition of the product in fact lies
in $D_{(n-k)-(-k)} = D_n$.  Similarly, note that $a^*$ is an element
of $L$ for each $a$ in $L$.

\statement{3.4. Proposition}{The product defined above is
associative.}

\proof If $a$ in is $D_n$ we shall denote by $a \delta_n$ the
element of $L$ given by \ $(a \delta_n) (m) = \delta_{n,m} a$, \ where
$\delta_{n,m}$ is the Kronecker symbol.

It is readily seen that the associativity of our product follows
from the identity
$$ (a_n \delta_n * a_m \delta_m) * a_p \delta_p = a_n \delta_n *
(a_m \delta_m * a_p \delta_p)$$
where $a_i \in D_i$ for $i = n, m, p$, which we now propose to
prove.  Using Proposition (3.3), the left hand side above becomes
$$ (\theta^n (\theta^{-n} (a_n) a_m) \delta_{n+m}) * a_p \delta_p =
$$
$$ \theta^{n+m} (\theta^{-n-m} (\theta^n (\theta^{-n} (a_n) a_m))
a_p) \delta_{n+m+p} =$$
$$ \theta^{n+m} (\theta^{-m} (\theta^{-n} (a_n) a_m) a_p)
\delta_{n+m+p}.$$

On the other hand, the right hand side of our identity equals
$$ a_n \delta_n * (\theta^m (\theta^{-m} (a_m) a_p) \delta_{m+p})=
\theta^n (\theta^{-n} (a_n) \theta^m (\theta^{-m} (a_m) a_p))
\delta_{n+m+p}.$$

Note that the term within the outermost parenthesis, to the right of
the last equal sign, is in $D_{-n} \cap D_m$, so that the coefficient
of $\delta_{n+m+p} $ above is in $D_{m+n} $ by (3.3).  Thus, proving
our identity, amounts to verify that
$$ \theta^{-m} (\theta^{-n} (a_n) a_m) a_p = \theta^{-n-m}
(\theta^n(\theta^{-n} (a_n) \theta^m (\theta^{-m} (a_m) a_p))$$
or, again by (3.3), that
$$ \theta^{-m} (\theta^{-n} (a_n) a_m) a_p = \theta^{-m}
(\theta^{-n} (a_n) \theta^m (\theta^{-m} (a_m) a_p)).$$

Let $(u_i)_{i}$ be an approximate identity for $D_{-m}$.  So the
left hand side above can be written as
$$ \lim_i \theta^{-m} (\theta^{-n} (a_n) a_m) u_i a_p = \lim_i
\theta^{-m} (\theta^{-n} (a_n) a_m \theta^m (u_i a_p)) =$$
$$ \lim_i \theta^{-m} (\theta^{-n} (a_n) \theta^m (\theta^{-m} (a_m)
u_i a_p)) = \theta^{-m} (\theta^{-n} (a_n) \theta^m (\theta^{-m} (a_m)
a_p)),$$
concluding the proof.
\endproof

\statement{3.5. Proposition}{For $a$ and $b$ in $L$ one has $(ab)^*
= b^* a^*$.}

\proof Arguing as in the beginning of the previous proof it is
enough to verify that
$$ (a_n \delta_n * a_m \delta_m )^* = (a_m \delta_m)^* * (a_n
\delta_n)^*$$
for $a_n \in D_n $ and $a_m \in D_m$.  The left hand side equals, by
definition and by (3.3)
$$ (\theta^n (\theta^{-n} (a_n) a_m) \delta_{n+m})^* = \theta^{-n-m}
(\theta^n (\theta^{-n} (a_n) a_m)^*) \delta_{-n-m} = $$
$$ \theta^{-m} (\theta^{-n} (a_n) a_m)^* \delta_{-n-m} = \theta^{-m}
((a_m)^* \theta^{-n} (a_n)^* ) \delta_{-n-m}.$$
On the other hand,
$$ (a_m \delta_m)^* * (a_n \delta_n)^* = \theta^{-m}(a_m)^*
\delta_{-m} * \theta^{-n} (a_n)^* \delta_{-n} =$$
$$ \theta^{-m} (\theta^m (\theta^{-m} (a_m)^*) \theta^{-n} (a_n)^*)
\delta_{-m-n} =$$
$$ \theta^{-m} ((a_m)^* \theta^{-n} (a_n)^* ) \delta_{-m-n}.$$

\endproof

\smallskip The main difficulties being overcome, we now have:

\statement{3.6. Theorem}{$L$ is an involutive Banach algebra with
the above defined multiplication, involution and norm.}

\statement{3.7. Definition}{The covariance algebra for the partial
automorphism $\Theta = (\theta, I, J)$ is the \cstar-algebra $C^* (A,
\Theta)$ obtained by taking the enveloping \cstar-algebra of $L$.}

\smallskip Elementary examples of this construction are standard crossed
products by the group of integers and the algebra of $n \times n$
complex matrices.  The latter is obtained as the covariance algebra
for the partial automorphism mentioned after (3.1).

A slightly more elaborate example is what one gets by taking $A =
c_0 (\N), I = A, J$ the ideal formed by sequences with a leading zero
and, finally, $\theta$ the forward shift.  The covariance algebra, in
this case, can be shown to be the algebra of compact operators on
$\ell_2 (\N)$.  If, instead, we took $A$ to be the unitization of
$c_0({\bf N}), I = A$, J the set of elements in $A$ with a leading
zero coordinate, and $\theta$ the forward shift, then the covariance
algebra becomes the Toeplitz algebra, that is, the \cstar-algebra
generated by the forward shift on $\ell_2 (\N)$.  The above statements
follow as easy corollaries of our structure Theorem (4.21), once one
considers the circle action given, in each case, by conjugation by
diag$(1, z, z^2, ...)$, for $z \in \S1$. With respect to the Toeplitz
algebra, see also (6.6).

The major example we would like to present is related to actions of
the circle group on \cstar-algebras.  As we shall see, any
\cstar-algebra having an action of $\S1$, under relatively mild
hypothesis on the action, is the covariance algebra for a certain
partial automorphism of the algebra of fixed points under $\S1$.  This
will be the subject of section (4).

\statement{3.8. Proposition}{Let $(e_i)_{i}$ be a (bounded)
approximate identity for $A$ then $(e_i \delta_0)_{i}$ is an
approximate identity for $L$ and hence also for $C^* (A, \Theta)$.}

\proof For $a_n \in D_n$ we have
$$ \lim_i (e_i \delta_0) * (a_n \delta_n) = \lim_i e_i a_n \delta_n
= a_n \delta_n.$$
By taking adjoints if follows that $\lim_i (a_n \delta_n)*(e_i
\delta_0) = a_n \delta_n$. Since $\sup_i \| e_i \delta_0 \| < \infty$
the above implies the conclusion.
\endproof

\statement{3.9. Proposition}{The map $E\colon  a \in L \rightarrow a (0)
\delta_0$ is a contractive positive conditional expectation
\cite{\RieffelInduced} from $L$ onto the subalgebra $A\delta_0$ of
$L$.}

\proof To prove positivity let $a \in L$.  Then
$$ (a*a^*)(0)=\sum^\infty_{k = -\infty} \theta^k (\theta^{-k} (a(k))
a^* (-k)) =$$
$$ \sum_k \theta^k (\theta^{-k} (a(k)) \theta^{-k} (a(k)^*)) =
\sum_k a(k) a(k)^* \geq 0.$$

The remaining statements can be easily verified and are left to the
reader.
\endproof

\statement{3.10. Corollary}{The obvious inclusion of $A$ into $L$,
composed with the map from $L$ into $C^*(A, \Theta)$, gives an
isometric *-homomorphism of $A$ into the latter.}

\proof Let $a\in A$ and let $f$ be a state on $A$ such that $f (a^*
a) = \| a \|^2$.  Identifying $A$ and its copy $A\delta_0$ within $L$
provides us with a state $f$ on $A\delta_0$ such that $f((a\delta_0)^*
* (a\delta_0)) = \| a \|^2$.  If that state is composed with the
conditional expectation of Proposition (3.9) we get a state on $L$.
By (3.8) $L$ has an approximate identity so one is allowed to use the
GNS construction, which provides a representation $\pi$ of $L$, and
hence of the covariance algebra, having a cyclic unit vector $\xi$ and
which satisfies
$$ \langle \pi (a \delta_0) \xi, \pi (a\delta_0) \xi \rangle = f
(a^*a) = \| a \|^2.$$

Thus, the norm $\| a \delta_0\|$, computed in the covariance
algebra, is no less than $\| a \|$.  The converse inequality follows
from the fact that the map mentioned in the statement is a
\cstar-algebra homomorphism, hence contractive.
\endproof

\smallskip We would now like to define the dual action, a concept
closely related to dual actions for crossed products.

Let, for every $z\in \S1$, $\alpha_z$ be the transformation of $L$
defined by
$$ (\alpha_z (a))(n) = z^n a(n) \for a\in L, z\in \S1.$$

The reader can easily verify that each $\alpha_z$ is a
*-automorphism of $L$ which, in turn, extends to a *-automorphism of
\cstar $(A, \Theta)$.  The resulting map $z \rightarrow \alpha_z$
becomes an action of $\S1$ on \cstar $(A, \Theta)$ which we shall call
the dual action.

\statement{3.11. Proposition}{For each $n$, let $B_n$ be the $n\th$
spectral subspace for the dual action.  Then the map $\phi_n\colon  x\in D_n
\rightarrow x\delta_n \in C^* (A, \Theta)$ is a linear isometry onto
$B_n$.}

\proof It is clear that $x \delta_n \in B_n$ for every $x \in D_n$.
Note that for such an $x$
$$ (x\delta_n) * (x \delta_n)^* = (x \delta_n) * (\theta^{-n} (x^*)
\delta_{-n} ) = \theta^n (\theta^{-n} (x) \theta^{-n} (x^*) ) \delta_0
= xx^* \delta_0.$$
Thus, in order to show that $\phi_n$ is an isometry, it suffices to
consider the case $n = 0$.  But this is just the conclusion of (3.10).
It now remains to show that the image of $\phi_n$ is all of $B_n$.  So
let $y \in B_n$ and write $y = \lim_k y_k$ where each $y_k$ belongs to
$L$ (or rather, the dense image of $L$ in $C^* (A, \Theta)$).  Note
that $L$ is invariant under the spectral projections of the covariance
algebra and also that $ y = P_n (y) = \lim_k P_n (y_k)$.  So, we may
assume that the $y_k$'s belong to the $n\th$ spectral subspace for the
corresponding action of $\S1$ on $L$.  That spectral subspace is
obviously $D_n \delta_n$ which, by our previous remarks, embeds
isometrically into the covariance algebra.  This implies that
$(y_k)_k$ is a Cauchy sequence with respect to the norm of $L$ and
hence that $y = \lim y_k \in D_n \delta_n$.
\endproof

\beginsection 4. The Structure of Actions of the Circle Group

In this section we intend to prove that any \cstar-algebra admitting
an action of $\S1$ is isomorphic to a covariance algebra, as described
above, provided the action satisfies two conditions which we now
describe.

\statement{4.1. Definition}{An action $\alpha$ of $\S1$ on a
\cstar-algebra $B$ is called semi-saturated if $B$ is generated, as a
\cstar-algebra, by the union of the fixed point algebra $B_0$ and the
first spectral subspace $B_1$ (compare \cite{\FellI}, \cite{\FellII},
\cite{\KishimotoTakai}, \cite{\PaschkeCircle}, \cite{\Phillips},
\cite{\RieffelProper}).}

\statement{4.2. Definition}{An action $\alpha$ of $\S1$ on a
\cstar-algebra $B$ is said to be stable if there exists an action
$\alpha^\prime$ on a \cstar-algebra $B^\prime$ such that $B \simeq
B^\prime \otimes K$ and $\alpha$ is the tensor product of
$\alpha^\prime$ by the trivial $\S1$-action on the algebra $K$ of
compact operators on a separable, infinite dimensional Hilbert space.}

\smallskip From now on we shall mainly be concerned with semi-saturated
stable actions on separable \cstar-algebras, but we would like to
argue that the above restriction are quite mild ones.  First of all
any action gives rise to a stable action by tensoring the old action
with the trivial action on $K$.  Moreover, if the old action was
semi-saturated then so will be the tensor product action.

Of course not all actions of $\S1$ are semi-saturated but this
difficulty could be circumvented in some cases.  That is, given an
action $\alpha$ on a \cstar-algebra $B$, it follows from (2.5) that
$B$ is generated by the union of all its spectral subspaces.

Now, if we let, for each positive integer $n$, $B^{(n)}$ be the
sub-\cstar-algebra of $B$ generated by $B_0$ and $B_n$ then $B^{(n)}$
is invariant under $\alpha$ and the formula
$$ \alpha^{(n)}_z = \alpha_{z^{1/n}} \for z \in S_1$$
provides a well defined semi-saturated action on $B^{(n)}$.

With some luck, one could put together the information obtained
about each $B^{(n)}$ to learn something of interest about the original
action.

\statement{4.3. Proposition}{Let $\alpha$ be a stable action of
$\S1$ on a separable \cstar-algebra $B$.  Then there exists an
isomorphism $\theta\colon  B_1^* B_1 \rightarrow B_1 B_1^*$ (see 2.2) and a
linear isometry $\lambda$ from $B^*_1$ onto $B_1 B_1^*$ such that for
$x, y \in B_1$, $a \in B_1^* B_1$ and $b \in B_1 B_1^*$
\medskip
\itemitem{(i)} $\lambda (x^* b) = \lambda(x^*) b$
\medskip
\itemitem{(ii)} $\lambda (ax^*) = \theta (a) \lambda (x^*)$
\medskip
\itemitem{(iii)} $\lambda (x^*)^* \lambda (y^*) = xy^*$
\medskip
\itemitem{(iv)} $\lambda (x^*) \lambda (y^{*})^{*} = \theta (x^*
y).$}

\proof Everything will follow from \cite{\BGR}, Theorem (3.4) after
we verify that both $B^*_1 B_1$ and $B_1 B^*_1$ are stable
\cstar-algebras with strictly positive elements.  Now, stability
follows at once from our assumption that the action be stable while
the existence of strictly positive elements is a consequence of our
separability hypothesis.
\endproof

\smallskip The conclusions of the last proposition will be among our
main tools in what follows.  In fact the only reason we shall consider
stable actions on separable \cstar-algebras is to obtain the
conclusions of (4.3).

If, for any other reason, the existence of $\theta$ and $\lambda$ as
above, are guaranteed in a specific example, we may discard stability
of the action $\alpha$ and separability of the algebra $B$ without
hurting the results of this section.

In fact it is quite common to find examples in which the conclusions
of (4.3) hold, but still the action is not stable.  The simplest such
example is provided by the action of $\S1$ on $B = M_2 ({\bf C})$
given by conjugation by the unitary matrices $\pmatrix{1 & 0 \cr 0 &
z}$ for $z$ in $\S1$.  In this case we have
$$ B_0 = \left\{\pmatrix{x & 0 \cr 0 & y} : x,y \in {\bf C}
\right\},$$
$$ B_1 = \left\{\pmatrix{0 & 0 \cr y & 0} : y \in {\bf C}
\right\},$$
$$ B_1^* B_1 = \left\{\pmatrix{x & 0 \cr 0 & 0} : x \in {\bf C}
\right\},$$
and
$$ B_1 B_1^* = \left\{\pmatrix{0 & 0 \cr 0 & y} : y \in {\bf C}
\right\}.$$

The maps $\theta$ and $\lambda$ defined by
$$ \theta \pmatrix{x & 0 \cr 0 & 0} = \pmatrix{0 & 0 \cr 0 & x} \for
x \in {\bf C}$$
and
$$ \lambda \pmatrix{0 & x \cr 0 & 0} = \pmatrix{0 & 0 \cr 0 & x}
\for x \in {\bf C}$$
satisfy the conclusions of (4.3), as the reader may easily verify,
although our action is obviously not stable.

\statement{4.4. Definition}{An action of $\S1$ on a \cstar-algebra
$B$ will be called regular if the conclusions of (4.3) hold. That is,
there should exist an isomorphism $\theta\colon  B_1^* B_1 \rightarrow B_1
B_1^*$ and a linear isometry $\lambda$ from $B^*_1$ onto $B_1 B_1^*$
satisfying (i) -- (iv) of (4.3).}

\smallskip We obviously have:

\statement{4.5. Corollary}{Stable circle actions on separable
\cstar-algebras are regular.}

\smallskip We shall also note the following fact.

\statement{4.6. Proposition}{Let $\alpha$ be the dual action on the
covariance algebra of a partial automorphism.  Then $\alpha$ is
regular.}

\proof The partial automorphism is built into the picture and
$\lambda$ is given by
$$ \lambda ((x \delta_1)^*) = x^* \delta_0 \for x \in D_1.$$

\endproof

\smallskip Also, there is no hope for an action of $\S1$ to be
equivalent to a dual action unless it is semi-saturated as we now
demonstrate.

\statement{4.7. Proposition}{If $\Theta = (\Theta, I, J)$ is a
partial automorphism of the \cstar-algebra $A$, then the dual action
$\alpha$ on the covariance algebra $B = C^* (A, \Theta)$ is
semi-saturated.}

\proof By (2.5) and by induction it is enough to show that, for $n
\geq 1$, $B_{n+1} \subseteq B_1 B_n$.  If $x \in B_{n+1}$ is of the
form $x = a \delta_{n+1}$ for $a\in D_{n+1}$, note that $\theta^{-1}
(a) \in D_{-1} \cap D_n$.  Without loss of generality we may assume
that $\theta^{-1} (a) = bc$ where $b, c \in D_{-1} \cap D_n$ (the
reduction to this case being made by an approximation argument).  Put
$b_1 = \theta(b)$ and observe that
$$ (b_1 \delta_1) * (c \delta_n) = \theta (\theta^{-1} (b_1) c)
\delta_{n+1} = \theta (bc) \delta_{n+1} = a \delta_{n+1} = x.$$
This shows that $B_{n+1}$, which equals $D_{n+1} \delta_{n+1} $ by
(3.11), is contained in $B_1 B_n$.
\endproof

\smallskip One of the crucial properties of semi-saturated actions is
obtained in our following result. It also provides an alternate
characterization of semi-saturation.

\statement{4.8. Proposition}{An action $\alpha$ of $\S1$ on $B$ is
semi-saturated if and only if for each $n > 0$ one has $B_n =
(B_1)^n$.  In this case, if $n, m \geq 0$ then $B_{n+m} = B_n B_m$.}

\proof Let $B^\prime $ be the dense *-subalgebra of $B$ generated by
$B_0 \cup B_1$.  Any element of $B^\prime $ is a sum of ``words'' of
the form $w = x_1 x_2 \ldots x_k$ where $x_i \in B_0 \cup B_1 \cup
B_1^*$.

A word, such as above, represents an element of $B$ which belongs to
some spectral subspace, the precise determination of which is obtained
by subtracting the number of $x_i$'s belonging to $B^*_1$ from the
number of those in $B_1$.  For this reason, to apply the spectral
projection $P_n$ to a sum of words is to eliminate all summands, but
those which belong to $B_n$.

If $b \in B_n$ let $b = \lim_j b_j$ with $b_j \in B^\prime $.  Thus
$b = P_n (b) = \lim_j P_n (b_j)$, so we may assume (upon replacing
$b_j$ by $ P_n (b_j))$ not only that $b_j \in B^\prime \cap B_n$, but
also that each $b_j$ is a sum of words each of which belongs to $B_n$.
A minute of reflection will show that all words in $B_n$ are in
$(B_1)^n$ for $n > 0$.  The converse follows from (2.5).
\endproof

\smallskip Let us now fix, for the remainder of this section, a
\cstar-algebra $B$ where a regular semi-saturated action $\alpha$ of
$\S1$ is defined.

\statement{4.9. Definition}{If $\lambda$ and $\theta$ are as in
(4.4) let
\ $\rho \colon  B^*_1 \rightarrow B_1^* B_1$,
\ $\lambda^\dagger\colon  B_1 \rightarrow B_1^* B_1$
\ and
\ $\rho^\dagger\colon  B_1 \rightarrow B_1 B^*_1$
\ be the isometries defined by
\ $\rho = \theta^{-1} \circ \lambda$,
\ $\lambda^\dagger (x) = \rho (x^*)^{*}$
\ and
\ $\rho^\dagger(x) = \lambda (x^*)^{*}$
\ for all $x \in B_1$.}

\statement{4.10. Proposition}{For $a \in B_1^* B_1, b \in B_1 B^*_1$
and $x, y \in B_1$ one has
\medskip
\itemitem{(i)} $ \rho(ax^*) = a\rho (x^*)$
\medskip
\itemitem{(ii)} $ \rho(x^*b) = \rho (x^*) \theta^{-1} (b)$
\medskip
\itemitem{(iii)} $ \rho(x^*) \rho (y^*)^{*} = x^* y$
\medskip
\itemitem{(iv)} $ \rho (x^*)^{*} \rho (y^*) = \theta^{-1} (x y^*)$
\bigskip
\itemitem{(v)} $ \lambda^\dagger (bx) = \theta^{-1} (b)
\lambda^\dagger (x)$
\medskip
\itemitem{(vi)} $ \lambda^\dagger (xa) = \lambda^\dagger (x) a$
\medskip
\itemitem{(vii)} $ \lambda^\dagger (x)^* \lambda^\dagger (y) = x^*
y$
\medskip
\itemitem{(viii)} $ \lambda^\dagger (x) \lambda^\dagger (y)^* =
\theta^{-1} (xy^*)$
\bigskip
\itemitem{(ix)} $ \rho^\dagger (bx) = b\rho^\dagger (x)$
\medskip
\itemitem{(x)} $ \rho^\dagger (xa) = \rho^\dagger (x) \theta (a)$
\medskip
\itemitem{(xi)} $ \rho^\dagger (x) \rho^\dagger (y)^* = xy^*$
\medskip
\itemitem{(xii)} $ \rho^\dagger (x)^* \rho^\dagger (y) = \theta (x^*
y).$}

\proof Follows by routine computations.
\endproof

\statement{4.11. Proposition}{The maps $\lambda, \rho,
\lambda^\dagger$ and $\rho^\dagger$ extend to isometries
$$ \lambda\colon  B^*_1 B \rightarrow B_1 B $$
$$ \rho\colon  B B^*_1 \rightarrow B B_1 $$
$$ \lambda^\dagger\colon  B_1 B \rightarrow B_1^* B $$
$$ \rho^\dagger\colon  B B_1 \rightarrow B B^*_1 $$
satisfying
\medskip
\itemitem{(i)} $\lambda (sb) = \lambda(s) b \for s \in B_1^* B_1,
b\in B$
\medskip
\itemitem{(ii)} $\lambda (s)^* \lambda (t) = s^* t \for s, t \in
B^*_1 B$
\medskip
\itemitem{(iii)} $\rho (bs) = b\rho(s) \for s \in B B^*_1 , b \in B$
\medskip
\itemitem{(iv)} $\rho (s) \rho (t)^* = st^* \for s, t \in B B^*_1$
\medskip
\itemitem{(v)} $\lambda^\dagger (sb) = \lambda^\dagger (s) b \for s
\in B_1 B, b \in B$
\medskip
\itemitem{(vi)} $\lambda^\dagger (s)^* \lambda^\dagger (t) = s^* t
\for s, t \in B_1 B$
\medskip
\itemitem{(vii)} $\rho^\dagger (bs) = b \rho^\dagger (s) \for s \in
B B_1, b \in B$
\medskip
\itemitem{(viii)} $\rho^\dagger (s) \rho^\dagger (t)^* = s t^* \for
s, t \in B B_1$.
\medskip
\noindent Moreover
\medskip
\itemitem{(a)} If $X \subseteq B_1^* B$ is a closed subspace,
invariant under left multiplication by $B^*_1 B_1$, then $\lambda (X)
= B_1 X$.
\medskip
\itemitem{(b)} If $X \subseteq B B^*_1$ is right $B_1
B^*_1$-invariant then $\rho (X) = XB_1$.
\medskip
\itemitem{(c)} If $X \subseteq B_1 B$ is left $B_1 B_1^*$-invariant
then $\lambda^\dagger (X) = B_1^* X$.
\medskip
\itemitem{(d)} If $X \subseteq B B_1$ is right $B^*_1 B_1$-invariant
then $\rho^\dagger (X) = XB^*_1$.
\medskip
\noindent Finally, $\rho^\dagger = \rho^{-1}$ and $\lambda^\dagger =
\lambda^{-1}$ hold for the extended maps.}

\proof If $s\in B^*_1 B$ is of the form $s = \sum^n_{i=1} x_i^* b_i$
where $x_i \in B_1$ and $b_i \in B$, put $\lambda (s) = \sum^n_{i=1}
\lambda (x_i^*) b_i$.  Observe that
$$ \| \sum^n_{i=1} \lambda (x_i^*) b_i \|^2 = \| \sum_{ij} b^*_j
\lambda (x_j^*)^* \lambda (x_i^*) b_i \| = \| \sum_{ij} b_j^* x_j
x_i^* b_i \| = \| \sum^n_{i=1} x_i^* b_i \|^2,$$
so that $\lambda$ is well defined, clearly satisfies (i) and (ii)
and is norm preserving.  The extensions of $\rho, \lambda^\dagger$ and
$\rho^\dagger$ are handled in a similar way.

Let $X$ be a closed subspace of $B^*_1 B$, invariant under left
multiplication by $B^*_1 B_1$.  Then by (2.6) $X = B^*_1 B_1 X$ so
$$ \lambda (X) = \lambda (B_1^*) B_1 X = B_1 B^*_1 B_1 X = B_1 X.$$

Similar proofs apply to the next three statements.

Finally, let $x, y \in B_1$ and note that
$$ \theta (y^*x) = \lambda (y^*) \lambda (x^*)^{*} = \lambda (y^*
\lambda (x^*)^{*} ) = \lambda (y^* \rho^\dagger (x)),$$
so that
$$ y^* x = \theta^{-1}(\lambda ( y^* \rho^\dagger (x))) = \rho (y^*
\rho^\dagger (x)) = y^* \rho (\rho^\dagger (x)).$$
Therefore, by (2.6), $x = \rho \rho^\dagger (x)$.

For $s = bx$ with $b\in B$ and $x \in B_1$ we have $\rho
\rho^\dagger (bx) = \rho (b \rho^\dagger (x)) = b \rho \rho^\dagger
(x) = bx $ so the composition
$$ B B_1 \stackrel{\rho^\dagger}{\rightarrow} B B_1^*
\stackrel{\rho}{\rightarrow} B B_1 $$
is the identity.  Since $\rho$ is isometric, hence injective, we get
$\rho^\dagger = \rho^{-1}$.  The proof that $\lambda^\dagger =
\lambda^{-1}$ goes along identical lines.
\endproof

\statement{4.12. Proposition}{The composition
$$ B^*_1 B B_1 \stackrel{\lambda}{\rightarrow} B_1 B B_1
\stackrel{\rho^\dagger}{\rightarrow} B_1 B B_1^*$$
is a \cstar-algebra isomorphism which extends $\theta$ and whose
inverse is given by the composition
$$ B_1 B B^*_1 \stackrel{\rho}{\rightarrow} B_1 B B_1
\stackrel{\lambda^\dagger}{\rightarrow} B_1^* B B_1 .$$}

\proof For $x, y \in B_1$ and $a \in B$ we have
$$ \rho^\dagger (\lambda (x^* ay)) = \rho^\dagger (\lambda (x^*) ay)
= \lambda (x^*) a \rho^\dagger (y).$$
Thus if $x_1, y_1 \in B_1$ and $a_1 \in B$
$$ \rho^\dagger (\lambda (x^* ay)) \rho^\dagger (\lambda(x^*_1 a_1
y_1)) = \lambda (x^*) a\rho^\dagger (y) \lambda (x^*_1) a_1
\rho^\dagger (y) = $$
$$ \lambda (x^*) a \lambda (y^*)^{*} \lambda (x_1^*) a_1
\rho^\dagger (y) = \lambda (x^*) a y x^*_1 a_1 \rho^\dagger (y) =
\rho^\dagger (\lambda ( x^* ay x_1^* a_1 y_1)), $$
proving $\rho^\dagger \circ \lambda$ to be multiplicative.

We also have
$$ \rho^\dagger (\lambda (x^* ay))^* = \rho^\dagger (y)^* a^*
\lambda (x^*)^{*} = \lambda (y^*) a^* \rho^\dagger (x) = \rho^\dagger
(\lambda (y^* a^* x)),$$
which shows that $\rho^\dagger \circ \lambda$ is star preserving.

To check the statement about extending $\theta$ note that for $x, y
\in B_1$ we have
$$ \rho^\dagger (\lambda (x^* y)) = \lambda (x^*) \rho^\dagger (y) =
\lambda (x^*) \lambda (y^*)^{*} = \theta (x^* y).$$
The fact that $\lambda^\dagger \rho$ inverts $\rho^\dagger \lambda$
follows from (4.11).
\endproof

\smallskip Our next result resembles one of the main axioms for
multipliers of \cstar-algebras.

\statement{4.13. Proposition}{If $s\in B B^*_1$ and $t \in B^*_1 B$
then $ \rho (s) t = s \lambda (t).$}

\proof It clearly suffices to verify the case in which $s = x^*$ and
$t = y^*$ with $x, y \in B_1$.  We have
$$ \rho (x^*) y^* = \rho (x^*) \rho (\rho^\dagger (y))^{*} = x^*
\rho^\dagger (y)^* = x^* \lambda (y^*).$$
\endproof

\statement{4.14. Proposition}{Regarding the partial automorphism
$(\theta, B_1^* B_1, B_1 B^*_1)$ of $B_0$ we have $\Dom (\theta^n) =
B^*_n B_n$ for all integers $n$.}

\proof This is clearly true for $n = -1, 0, 1$.  For $n \geq 1$,
using induction, we have
$$ \Dom (\theta^{n+1}) = \{x \in \Dom (\theta): \theta(x) \in \Dom
(\theta^n)\} = \{x \in B_1^* B_1: \theta(x) \in B^*_n B_n \}.$$
Thus, given $x$ in $\Dom (\theta^{n+1})$
$$ \theta(x) \in \Im (\theta) \cap \Dom (\theta^n) = \Im (\theta) \cdot
\Dom (\theta^n) = B_1 B_1^* B_n^* B_n B_1 B_1^*,$$
so
$$ x = \theta^{-1} (\theta (x)) = \lambda^\dagger \rho (\theta(x))
\in \lambda^\dagger (\rho (B_1 B_1^* B_n^* B_n B_1 B_1^*)) =$$
$$ B_1^* B_1 B_1^* B_n^* B_n B_1 B_1^* B_1 = B_{n+1}^* B_{n+1}.$$
This shows that $\Dom (\theta^{n+1}) \subseteq B^*_{n+1} B_{n+1}$.
Conversely, if $x \in B_{n+1}^* B_{n+1}$ then clearly $x \in \Dom
(\theta) $ and
$$ \theta (x) = \rho^\dagger \lambda (x) \in \rho^\dagger (\lambda
(B_{n+1}^* B_{n+1})) = B_1 B_{n+1}^* B_{n+1} B^*_1 \subseteq B^*_n B_n
= \Dom (\theta^n).$$
So $x \in \Dom (\theta^{n+1})$.
\endproof

\smallskip Our main goal, as we already indicated, is to describe $B$ as
the covariance algebra for the partial automorphism $\Theta = (\theta,
B^*_1 B_1, B_1 B^*_1)$ of the fixed point algebra $B_0$.  The
conclusions of (4.14) suggest that, if we are to succeed, then $B_n
B_n^*$, being essentially the $n\th$ spectral subspace of the
covariance algebra $C^* (B_0, \Theta)$, should correspond to $B_n$.
This leads us to the following definitions.

\statement{4.15. Definition}{Let $i_0\colon  B_0 \rightarrow B_0$ denote
the identity map of $B_0$ and
$i_n\colon  B_n \rightarrow B_0$
be defined, inductively, by
$i_n (x) = i_{n-1} (\rho^\dagger (x))$
when $n \geq 1$ or by
$i_n (x) = i_{n+1} (\rho (x))$
when $n \leq -1$.}

\smallskip It is clear that $i_n$ is an isometry for each $n$ (although
not necessarily surjective).

\statement{4.16. Proposition}{For all integers $m$ and $n$ and for
any $x_n \in B_n$ and $y_m \in B_m$ we have
$$ i_{n+m} (x_n y_m) = i_n (x_n i_m (y_m)).$$}

\proof Assume initially that $n = 1$.  If $m=0$ there is nothing to
prove and if $m > 0$, using induction, we have
$$ i_{1+m} (x_1 y_m) = i_m (\rho^\dagger (x_1 y_m)) = i_m (x_1
\rho^\dagger (y_m)) = $$
$$ i_1 (x_1 i_{m-1} \rho^\dagger (y_m)) = i_1 (x_1 i_m (y_m)).$$

If $m = -1$, still assuming $n = 1$, we have
$$ i_1 (x_1 i_{-1} (y_{-1})) =\rho^\dagger (x_1 \rho (y_{-1})) =
\rho^\dagger \rho (x_1 y_{-1}) = x_1 y_{-1} = i_0 (x_1 y_{-1}).$$

Now let $ m < -1$ and observe that, by induction
$$ i_{1+m} (x_1 y_m) = i_{1+(m+1)} (\rho ( x_1 y_m)) = i_{1+(m+1)}
(x_1 \rho (y_m)) =$$
$$ i_1 (x_1 i_{m+1} \rho (y_m)) = i_1 (x_1 i_m (y_m)).$$

If $n > 1$ we may assume, without loss of generality, by (4.8) or
(2.6) that $x_n = x_1 x_{n-1}$ where $x_1 \in B_1$ and $x_{n-1} \in
B_{n-1}$.  So, by induction once more,
$$ i_{n+m} (x_n y_n) = i_{1+(n-1)+m} (x_1 x_{n-1} y_m) = i_1 (x_1
i_{n-1+m} (x_{n-1} y_m)) =$$
$$ i_1 (x_1 i_{n-1} (x_{n-1} i_m (y_m))) = i_n (x_1 x_{n-1} i_m
(y_m)).$$

The case $n = 0$ is trivial and the proof for $n < 0$ is done in a
similar way.
\endproof

\statement{4.17. Proposition}{For every integer $n$ one has
$i_n(B_n) = B_n B^*_n$.}

\proof The statement follows by definition when $n = -1, 0, 1$.  If
$n \geq 1$, by induction, (2.7), (4.8) and (4.11) one has
$$ i_{n+1} (B_{n+1}) = i_1 ( B_1 i_n (B_n)) = \rho^\dagger (B_1 B_n
B_n^*)= $$
$$ \rho^\dagger (B_1 B^*_1 B_1 B_n B_n^*) = \rho^\dagger (B_1 B_n
B_n^* B_1^* B_1) = B_1 B_n B_n^* B^*_1 B_1 B_1^* = B_{n+1}
B_{n+1}^*.$$

As usual, the proof for negative values of $n$ is omitted.
\endproof

\statement{4.18. Proposition}{If $n\in \Z$, $a\in B^*_n B_n$ and
$x_n \in B_n$ then $i_n(x_n a) = i_n(x_n) \theta^n(a).$}

\proof For $n=-1,0,1$ the statement is obvious. For $n>1$ assume, as
before, that $x_n = x_1 x_{n-1}$, the subscripts indicating the
spectral subspace each factor lies in. So, by induction,
$$ i_n(x_na) = i_n(x_1 x_{n-1} a) = i_1(x_1 i_{n-1}(x_{n-1}a)) =$$
$$ \rho^\dagger(x_1 i_{n-1}(x_{n-1})\theta^{n-1}(a)) = \rho^\dagger
(x_1 i_{n-1}(x_{n-1}))\theta^n(a) = i_n(x_n)\theta^n(a).$$

For $n<-1$ the proof is similar.
\endproof

\statement{4.19. Proposition}{For all integers $n$ and $m$ and for
all $x_n \in B_n$ and $x_m \in B_m$}
$$ \theta^{-n} (i_{n+m} (x_n y_m)) = \theta^{-n} (i_n (x_n)) i_m
(x_m).$$

\proof Let $\{e_i \}$ be an approximate identity for $B_n^* B_n$.
We have
$$ \theta^{-n} (i_{n+m} (x_n y_m)) = \theta^{-n} (i_n(x_n i_m
(y_m))) = \lim_j \theta^{-n} (i_n (x_n e_j i_m (y_m))) =$$
$$ \lim_j \theta^{-n} (i_n (x_n) \theta^n (e_j i_m (y_m))) = \lim_j
\theta^{-n} (i_n (x_n)) e_j i_m (y_m)) =\theta^{-n} (i_n (x_n)) i_m
(y_m).$$

\endproof

\statement{4.20. Proposition}{For any $n$ and for any $x_n \in B_n$
one has $\theta^{-n} (i_n (x_n)) = i_{-n} (x_n^*)^*$.}

\proof This is obvious for $n = 0$.  If $n = 1$ we have
$$ \theta^{-1} (i_1 (x_1)) = \lambda^\dagger \rho \rho^\dagger (x_1)
= \lambda^\dagger (x_1) = \rho (x_1^*)^{*} = i_{-1} (x_1^*)^{*}.$$
If $ n \geq 1$ assume, without loss of generality, that $x_n = x_1
x_{n-1}$ where $x_1 \in B_1$ and $x_{n-1} \in B_{n-1}$.  So, by
induction, we have
$$ \theta^{-n} (i_n (x_n)) =\theta^{-(n-1)} \theta^{-1} (i_1 (x_1
i_{n-1} (x_{n-1}))) = \theta^{-(n-1)} (i_{-1} (i_{n-1}(x_{n-1})^*
x_1^*))^* =$$
$$ \theta^{-(n-1)} (i_{n-1} (x_{n-1})^* i_{-1} (x_1^*))^* =
\theta^{-(n-1)} (i_{-1} (x_1^*) i_{n-1} (x_{n-1})) =$$
$$ \theta^{-(n-1)} (i_{n-1} (i_{-1} (x_1^*) x_{n-1})) = i_{-(n-1)}
(x^*_{n-1} i_{-1} (x_1^*)^{*})^* =$$
$$ i_{-n} (x_{n-1}^* x^*_1)^{*} = i_{-n} (x_n^*)^{*}.$$

The proof for negative values of $n$ is similar.
\endproof

\smallskip We are now prepared to present our first main result.

\statement{4.21. Theorem}{Let $\alpha$ be a semi-saturated regular
action of $\S1$ on a \cstar-algebra $B$.  If $\Theta = (\theta, B_1^*
B_1, B_1 B^*_1)$ is the partial automorphism of the fixed point
algebra $B_0$ as in (4.4), then there exists an isomorphism
$$ \phi\colon  C^* (B_0, \Theta) \rightarrow B$$
which is covariant with respect to the dual action.}

\proof Recall that $L$ is the Banach *-algebra formed by summable
sequences $(a(n))_{n\in {\bf Z}}$ with $a(n) \in \Dom(\theta^{-n}) =
B_n B_n^*.$ Define
$$ \phi \colon a \in L \rightarrow \sum i^{-1}_n (a(n)) \in B.$$
To show that $\phi $ is multiplicative, let $a_n \in B_n B_n^*$,
$b_m \in B_m B_m^*$ and put $x_n = i_n^{-1} (a_n)$ and $y_m = i_m^{-1}
(b_m)$. We must then verify that
$$ \phi (a_n \delta_n * b_m \delta_m) = \phi (a_n \delta_n) \phi
(b_m \delta_m).$$
For this end, we need to check that
$$ i_{n+m}^{-1} (\theta^n (\theta^{-n} (a_n) b_m)) = i_n^{-1} (a_n)
i_m^{-1} (b_m)$$
or that
$$ \theta^n (\theta^{-n} (i_n (x_n )) i_m (y_m)) = i_{n+m} (x_n
y_m)$$
which follows from (4.19).

To check that $\phi$ preserves adjoints it is enough to show that,
for any $n$ and any $a\in B_n B_n^*$,
$$ \phi ((a \delta_n )^*) = \phi (a \delta_n)^*$$
which translates to
$$ i^{-1}_{-n} (\theta^{-n} (a^*)) = i^{-1}_n (a)^*.$$

Writing $a = i_n (x_n)$, for $x_n \in B_n$, the above becomes
$$ i^{-1}_{-n} (\theta^{-n} (i_n (x_n)^*)) = x^*_n$$
which is equivalent to the conclusion of (4.20).

If follows that $\phi$ extends to a *-homomorphism from $C^* (B_0,
\Theta)$ to $B$ which is surjective, because its image contains any
$B_n$, and which is equivariant with respect to the dual action, since
spectral subspaces are mapped accordingly.

It remains to show that $\phi$ is injective, but in view of (2.9) we
only need to check injectivity on the fixed point subalgebra of
$C^*(B_0, \Theta)$, a fact that follows from (3.11). \endproof

\beginsection 5. Representations of Covariance Algebras

The parallel between our covariance algebras and crossed product
algebras is reflected, also, in their representation theory which we
shall now study in detail.  In particular, it is possible to define an
analog for the regular representation of crossed products which is our
next goal.

Fix, for the remainder of this section, a \cstar-algebra $A$ and a
partial automorphism $\Theta = (\theta, I, J)$ of $A$.

\statement{5.1. Definition}{Let $\pi $ be a representation of $A$ on
the Hilbert space $\H$.  The regular representation $\tilde{\pi} $ of
$C^*(A, \Theta)$, associated to $\pi$, is that which is obtained by
inducing $\pi$ to $L$ \cite{\RieffelInduced} via the conditional
expectation $E$ of (3.9) and then extending it to the enveloping
\cstar-algebra of $L$.}

\smallskip A result which could perhaps be thought of as a new
manifestation of the amenability of the group of integers is the
following (see \cite{\Pedersen}):

\statement{5.2. Theorem}{If $\pi$ is a faithful representation of
$A$ then $\tilde{\pi} $ is faithful on $C^*(A, \Theta)$.}

\proof The inducing process \cite{\RieffelInduced} leads us to
consider the Hilbert space $\tilde{\H}$ obtained by completing $L
\otimes \H$ under the semi-norm given by the (sometimes degenerated)
inner product
$$ \langle y_1 \otimes \xi_1, y_2 \otimes \xi_2 \rangle = \langle
\pi (E(y^*_2 y_1)) \xi_1, \xi_2 \rangle \for y_i \in L, \xi_i \in \H,
i = 1, 2.$$
The induced representation itself is specified by
$$ \tilde{\pi} (x) (y \otimes \xi) = (xy) \otimes \xi \for x, y \in
L, \xi \in \H.$$

Let $L_n$ be the $n\th$ spectral subspace for the dual action
restricted to $L$, that is, $L_n = D_n \delta_n$.  Since, for $n \not
= m$, $E (L_n^* L_m) = 0$ and because $\bigoplus_n L_n$ is dense in
$L$, we have that the subspaces $\tilde{\H}_n$ of $\tilde{\H}$,
obtained by closing the image of $L_n \otimes \H$ in $\tilde{\H}$,
form an orthogonal decomposition of $\tilde{\H}$.

If $y = d_n \delta_n$ for $d_n \in D_n$, then it is clear that
$\tilde{\pi} (y)$ sends $\tilde{\H}_m$ into $\tilde{\H}_{m+n}$ for
every $m$.  Thus, if $y$ is the finite sum $y = \sum_{n=-N}^N d_n
\delta_n$, it follows that
$$ p_{m+n} \circ \tilde{\pi} (y) |_{\H_{m}} = \tilde{\pi} (P_n (y))
|_{\H_{m}} \for n, m \in \Z,$$
where $p_n$ is the orthogonal projection of $\tilde{\H}$ onto
$\tilde{\H}_n$ and $P_n$ is, as in (2.4), the $n\th$ spectral
projection.  Since the finite sums are dense in $C^*(A, \Theta)$, a
continuity argument yields the above formula for all $y \in C^*(A,
\Theta)$.

{}From this it can easily be seen that if $\tilde{\pi} (y) = 0$ then
$\tilde{\pi} (P_n (y)) = 0$ for any $y$ in the covariance algebra.
This reduces our task to the easier one of proving $\tilde{\pi}$ to be
injective on each spectral subspace of $C^*(A, \Theta)$, which we have
already shown to be essentially $D_n \delta_n$.  So let $d_n \in D_n$
and suppose that $\tilde{\pi}(d_n \delta_n) = 0$.  Note that for $a\in
D_0 = A$ and $\xi \in \H$
$$ 0 =\tilde{\pi} (d_n \delta_n) (a \delta_0 \otimes \xi) = \theta^n
(\theta^{-n} (d_n) a) \delta_n \otimes \xi.$$
Therefore, the right hand side represents the zero vector in
$\tilde{\H}$.  By definition, its norm in $\tilde{\H}$ is computed by
the horrible looking expression
$$ \langle \pi (E((\theta^n (\theta^{-n} (d_n)a) \delta_n)^* *
(\theta^n(\theta^{-n}(d_n)a)\delta_n))) \xi, \xi \rangle^{1/2}$$
which, fortunately, can be simplified to $\| \pi (\theta^{-n} (d_n)
a) \xi \|$.

Since $\xi$ can be any vector of $\H$, it follows that $\pi
(\theta^{-n} (d_n) a) = 0 $.  But $\pi$ was supposed faithful.
Letting $a$ run through an approximate unit for $A$ we get
$\theta^{-n} (d_n) = 0$ and hence $d_n = 0$.
\endproof

\smallskip Representations of crossed product algebras are known to
correspond to covariant representations of the corresponding
\cstar-dynamical system \cite{\Pedersen}.  With the necessary
modifications, the same is true for covariance algebras of partial
automorphisms as we shall now see.

\statement{5.3. Definition}{A covariant representation of the pair
$(A, \Theta)$ is a triple $(\H, \pi, u)$ where $\pi$ is a
*-representation of $A$ on the Hilbert space $\H$ and $u \in B(\H)$ is
a partial isometry whose initial space is $\pi (I)\H$ (meaning closed
linear span) and whose final space if $\pi (J)\H$.  In addition it is
required that, for $a \in I$,
$$ \pi (\theta(a)) = u \pi (a) u^*.$$}

\smallskip The rules of the game are often so tough when dealing with
algebraic properties of partial isometries that the reader may be
stricken by the naiveness of the above definition.  In fact that was
also the author's first impression. It is really quite surprising that
this definition happens to do it's job in a satisfactory way.

A covariant representation gives rise to a representation of $C^*(A,
\Theta)$ as we shall see. Before that, let us introduce some notation.

\statement{5.4. Definition} {If $n$ is a negative integer and if
$v$ is a partial isometry on a Hilbert space, then $v^n$ stands for
$v^{{-n}^*}$ (note that this will not cause any confusion since, when
$v$ is invertible, $v^* = v^{-1}$ and the two possible interpretations
of $v^n$ coincide).}

\smallskip For $y \in L $ (cf. 3.6) put
$$ (\pi\times u) (y) = \sum_{n\in {\bf Z}} \pi (y(n)) u^n.$$

\statement{5.5. Proposition}{$\pi\times u$ is a representation of
$L$ on $\H$ and therefore, extends to a representation of $C^*(A,
\Theta)$.}

\proof The proof of the multiplicativity of $\pi\times u$ consists
in justifying all steps of the following calculation, for $a_n \in
D_n$ and $a_m \in D_m$
$$ \pi (a_n) u^n \pi (a_m) u^m = u^n \pi (\theta^{-n}
(a_n) a_m) u^m = $$
$$ \pi (\theta^n (\theta^{-n} (a_n ) a_m)) u^n u^m = \pi
(\theta^n (\theta^{-n} (a_n) a_m )) u^{n+m}.$$
The first two steps follow from the formula
$$ u^n \pi (\theta^{-n} (a_n)) = \pi (a_n) u^n \for a_n \in
D_n,$$
which can be proved, using induction, separately for $n > 0$ and $n
< 0$.

The final step is a consequence of the formula,
$$ \pi (a_n) u^n u^m = \pi (a_n) u^{n+m} \for a_n \in
D_n,$$
which is obviously true when $n$ and $m$ have the same sign, but
which must be proved by induction otherwise.  The proof that
$\pi\times u$ is invariant under the star operation is left to the
reader.
\endproof

\statement{5.6. Theorem}{Let $\sigma$ be a *-representation of $C^*
(A, \Theta)$ on the Hilbert space $\H$.  Then there exists a covariant
representation $(\pi, \H, u)$ of the pair $(A, \Theta)$ such that
$\sigma = \pi\times u$.}

\proof Let $B = C^* (A, \Theta)$ and let $\lambda, \rho,
\lambda^\dagger $ and $\rho^\dagger$ be as in (4.11).  Denote by $\pi$
the restriction of $\sigma$ to $A$, identified with $A\delta_0$.  For
$x_i \in I$ and $\xi_i \in \H$, $i = 1, \ldots, n$ let
$$ u (\sum^n_{i=1} \pi (x_i) \xi_i) = \sum^n_{i=1} \sigma (\lambda
(x_i \delta_0)) \xi_i,$$
and let us verify that this extends to a well defined isometry from
$\pi (I)\H$ to $\pi (J) \H$.  We have
$$ \| \sum_i \sigma (\lambda (x_i \delta_0)) \xi_i \|^2 = \sum_{ij}
\langle \sigma ( \lambda (x_j \delta_0)^* * \lambda (x_i \delta_0) )
\xi_i, \xi_j \rangle = $$
$$ \sum_{ij} \langle \pi (x_j^* x_i) \xi_i, \xi_j\rangle = \| \sum_i
\pi (x_i) \xi_i \|,$$
so $u$ is well defined and isometric, but we still have to check
that the image of $u$ is $\pi (J) \H$.  That image equals
$$ \sigma (\lambda (B_1^* B_1)) \H = \sigma (B_1 B_1^* B_1) \H
\subseteq \pi (B_1 B_1^*) \H = \pi (J) \H,$$
while
$$ \pi (J) \H = \pi (B_1 B_1^*) \H = \pi (B_1 B^*_1 B_1 B_1^*)\H
\subseteq \sigma (B_1 B_1^* B_1) \H = \sigma (\lambda ( B_1^* B_1 ))
\H = \Im (u).$$

We next need to show that $(\pi, \H, u)$ satisfies the covariance
equation $\pi ( \theta (a)) = u \pi (a) u^*$ for $a \in I$.  Given $b
\in I$ and $\xi \in \H$ we have
$$ u \pi (a) \pi (b) \xi = u \pi (ab) \xi = \sigma (\lambda (ab
\delta_0 )) \xi = \sigma ((\theta (a) \delta_0) * \lambda (b
\delta_0)) \xi = \pi (\theta (a)) u (\pi (b) \xi),$$
which says that $u \pi (a) = \pi (\theta (a)) u$ on $ \pi (I) \H$.
But both $u \pi (a)$ and $\pi (\theta (a)) u$ vanish on the orthogonal
complement of $\pi (I)\H$ hence the equality above holds on the whole
of $\H$.  So
$$ u \pi (a) u^* = \pi (\theta (a)) uu^* = \pi (\theta(a)).$$

Recall from (4.6) that $\lambda ((a_1 \delta_1)^*) = a^*_1 \delta_0$
for $a_1 \in D_1$.  Given $b_{-1} \in D_{-1} = I$, of the form $b_{-1}
= \theta^{-1} (c_1^* d_1)$ with $c_1, d_1$ in $D_1$ note that
$$ (c_1 \delta_1)^* * (d_1 \delta_1) = \theta^{-1} (c_1^* d_1)
\delta_0 = b_{-1} \delta_0.$$
So
$$ \lambda (b_{-1} \delta_0) = \lambda ((c_1 \delta_1)^*) *(d_1
\delta_1) = c_1^* d_1 \delta_1 = \theta (b_{-1}) \delta_1.$$

The above form of $b_{-1} $ is sufficiently general to insure that
$$ \lambda (b_{-1} \delta_0) = \theta (b_{-1}) \delta_1$$
for any $b_{-1} \in D_{-1}$.  This equality is the basis of the
proof that $\pi\times u = \sigma$.  In fact, since the dual action is
semi-saturated as proven in (4.7), it is enough to show that
$\pi\times u$ and $\sigma$ coincide on $B_0 \cup B_1$.  This is
obvious for $B_0$ so it remains to prove that
$$ \sigma (b_1 \delta_1) = \pi (b_1) u \for b_1 \in D_1.$$

Both sides of the equal sign above represent operators which vanish
on the orthogonal complement of $\pi (I) \H$ so, in order to prove
their equality, we may restrict ourselves to $\pi (I) \H$.  Given
$a_{-1} \in I$ and $\xi \in \H$ we have for any $b_1 \in D_1$
$$ \pi (b_1) u \pi (a_{-1}) \xi = \sigma (b_1 \delta_0) \sigma(
\lambda (a_{-1} \delta_0 )) \xi = \sigma ((b_1 \delta_0)*( \theta
(a_{-1}) \delta_1) ) \xi =$$
$$ \sigma (\theta (\theta^{-1} (b_1) a_{-1} ) \delta_1 ) \xi =
\sigma ((b_1 \delta_1)*( a_{-1} \delta_0)) \xi = \sigma (b_1 \delta_1)
\pi (a_{-1}) \xi$$
which concludes our proof.
\endproof

\beginsection 6. The Toeplitz Algebra of a Partial Automorphism

Let $A$ be a \cstar-algebra and $\Theta = (\theta, I,J)$ be a partial
automorphism of $A$.  We denote by $B$, the covariance algebra $C^*
(A, \Theta)$ and regard $A$ as a subalgebra of $B$ in the obvious way.

Fixing a faithful representation of $B$, allows us to assume that
$B \subseteq B(\H)$ for some Hilbert space $\H$.  Using Theorem (5.6)
we conclude that there is a partial isometry $u$ in $B(\H)$ such that,
for $x \in I$ and $y \in J$, we have
\medskip
\itemitem{(i)} $u x u^* = \theta (x)$
\medskip
\itemitem{(ii)}$x u^* u = x$
\medskip
\itemitem{(iii)} $yuu^* = y$.
\medskip
\noindent Moreover, the set of finite sums of the form $\sum_{n \in
\Z} a_n u^n$ with $a_n \in D_n = \Dom(\theta^{-n})$, is dense in $B$
(see 5.4). The dual action of $\S1$ on $B$ is given, on the above
mentioned dense set, by
$$\alpha_z (\sum a_n u^n) = \sum z^n a_n u^n \for z \in \S1.$$
So the $n\th$ spectral subspace for $\alpha$ is given by $B_n = D_n
u^n$ (see 3.11).

In the following, $\N^*$ denotes the set of strictly positive
integers, $S$ denotes the unilateral shift on $\ell_2 (\N^*)$ while $P
= 1 - SS^*$ and $Q = SS^*$.  Also let $e_{ij}$ denote the standard
matrix units in $B(\ell_2 (\N^*))$ for $i, j \geq 1$.  Note that
$e_{ij} = S^i S^{*j} - S^{i-1} S^{*j-1}$.

\statement{6.1. Definition} {The Toeplitz algebra for the pair $(A,
\Theta)$ is the sub-\cstar-algebra $\E = \E (A, \Theta)$, of operators
on $\H\otimes \ell_2 (\N^*)$, generated by the set
$$\{ b_n \otimes S^n: n\in \Z, b_n \in B_n\}.$$}

\smallskip Since the dual action $\alpha$ is semi-saturated, we have by
(4.7) and (4.8), that $B_n = (B_1)^n$ for $n \geq 1$.  Our next result
is an immediate consequence of this.

\statement{6.2. Proposition} {$\E$ is generated by $(B_0 \otimes 1)
\cup (B_1 \otimes S)$.}

\proof It is enough to note that, for $n \geq 1$
$$(B_n \otimes S^n) = ((B_1)^n \otimes S^n) = (B_1 \otimes S)^n.$$
\endproof

\statement{6.3. Definition} {Let $\Lambda$ be the subset of
$B(\H\otimes \ell_2 (\N^*))$ given by
$$\Lambda = \bigoplus_{i,j\geq 1} (B_i B^*_j \otimes e_{ij}).$$}

\smallskip From now on we will use the symbol $\bigoplus$ to denote the
\underbar{closure} of the sum of any family of independent subspaces
of a Banach space.  Also, following (2.2), we regard $B_i B_j^*$ as
the closed linear span of the set of products $xy$ where $x \in B_i$
and $y \in B^*_j$.

Because, for any $j$, $B_j^* B_j \subseteq B_0$ and also because
$B_i B_0 \subset B_i$, we have that $\Lambda $ is a sub-\cstar-algebra
of $B(\H\otimes \ell_2 (\N^*))$.

\statement{6.4. Proposition} {$\Lambda$ is an ideal in $\E$.}

\proof Since $B_n = (B_1)^n$, for $n \geq 1$, the set of elements of
the form $b_{n-1} b_1 b_m^*$, the indices indicating the spectral
subspace of $B$ each factor lies in, is total (their linear span is
dense) in $B_n B_m^*$.  So the identity
$$b_{n-1} b_1 b_m^* \otimes e_{nm} = (b_{n-1} \otimes S^{n-1} ) (b_m
b_1^* \otimes S^{m-1})^*- (b_{n-1}b_1 \otimes S^n) (b_m \otimes
S^m)^*$$
implies that $\Lambda$ is contained in $\E$.  The verification that
$\Lambda$ is actually an ideal is now straightforward.  \endproof

\smallskip As in \cite{\PV}, one should now show that the quotient
$\E/\Lambda$ is isomorphic to $B$, providing the short exact sequence
$$0\rightarrow \Lambda \rightarrow \E \rightarrow B \rightarrow 0$$

Even though we could attempt to prove this right now, we shall
deduce it from a much stronger result (Lemma 6.7 below), which we will
need for various other reasons.

For each $z\in \S1$ let $\Delta_z$ be the infinite matrix
$$\Delta_z = \pmatrix{
  1 & & & \cr
  & z & & \cr
  & & z^2 \cr
  & & & \ddots
  }$$
regarded as a unitary operator on $\H\otimes \ell_2 (\N^*)$.  Note
that $\Delta_z (b_n \otimes S^n) \Delta_z^{-1} = z^n (b_n \otimes
S^n)$ so that $\Delta_z \E \Delta^{-1}_z = \E$ and so we can define an
action $\gamma$ of $\S1$ on $\E$ by
$$\gamma_z (x) = \Delta_z x \Delta_z^{-1} \for x\in \E, z \in \S1.$$

Our next goal is to show that $\gamma$ is a semi-saturated regular
action, which will enable us to describe $\E$ as the covariance
algebra for a certain partial automorphism, by Theorem (4.21).

\statement{6.5. Proposition} {Let $\E_n$ denote the $n\th$ spectral
subspace relative to $\gamma$.  Then
\medskip
\itemitem{(i)} $\E_0 = A \otimes 1 \ +\ \bigoplus^\infty_{n=1} (D_n
\otimes e_{nn})$
\medskip
\itemitem{(ii)} $\E_1 = B_1 \otimes S \ +\ \bigoplus^\infty_{n=1}
(B_{n+1} B_n^* \otimes e_{n+1, n})$
\medskip
\itemitem{(iii)} $\E^*_1 \E_1 = D_{-1} \otimes 1\ +\
\bigoplus^\infty_{n=1} ((D_{-1} \cap D_n) \otimes e_{nn})$
\medskip
\itemitem{(iv)} $\E_1 \E_1^* = D_1 \otimes Q \ +\
\bigoplus^\infty_{n=2} (D_n \otimes e_{nn})$.}

\proof Let $Y_0 $ and $Y_1$ denote the right hand side of (i) and
(ii), respectively.  Define $Y_n = (Y_1)^n$ for $n > 1$, and $Y_n =
(Y_1^*)^{-n}$ for $n < 0$.  It can be shown without much difficulty,
that $Y = \bigoplus_{n\in \Z} Y_n$ is a sub-\cstar-algebra of $\E$.

Should the reader decide to verify this by himself, we suggest he
starts by proving that $Y_1 Y_0$ and $Y_0 Y_1$ are contained in $Y_1$
and that $Y_1^* Y_1$ and $Y_1 Y_1^*$ are contained in $Y_0$.

We next note that $B_0 \otimes 1 \subseteq Y_0$ and that $B_1
\otimes S \subseteq Y_1$ so, by (6.2) we get $Y = \E$.  Since it is
obvious that $Y_n \subseteq \E_n$ we must have $Y_n = \E_n$ for all
$n$ and in particular for $n = 0$ and $n=1$, proving (i) and (ii).
The two remaining statements follow easily from (ii).  \endproof

\statement{6.6. Proposition} {$\gamma$ is a semi-saturated regular
action of $\S1$ on $\E$. The roles of the maps $\theta$ and $\lambda$
mentioned in (4.4) are played, respectively, by $\theta_\E$ and
$\lambda_\E$ defined by
$$\theta_\E \colon  x \in \E_1^* \E_1 \longmapsto (u \otimes S) x
(u\otimes S)^* \in \E_1 \E^*_1$$
$$\lambda_\E \colon  x^* \in \E_1^* \longmapsto (u \otimes S) x^* \in \E_1
\E^*_1.$$
Therefore $\E$ is isomorphic to the covariance algebra
$C^*(\E_0,(\theta_\E, \E_1^*\E_1, \E_1 \E_1^*)).$}

\proof The proof of regularity consists in verifying that the above,
in fact, describes well defined maps between the indicated sets and
that they satisfy (4.4).  These are somewhat routine tasks, which we
leave for the reader. Finally, since it is obvious that $B_0 \otimes 1
\subseteq \E_0$ and that $B_1 \otimes S \subseteq \E_1$,\ $\gamma$ is
seen to be semi-saturated by (6.2).  \endproof

\smallskip Our next Lemma is an important tool in obtaining
representations of the Toeplitz algebra (compare Lemma 7.1 in
\cite{\AP}).

\statement{6.7. Lemma} {Let $\K$ be a Hilbert space and $\pi$ be a
representation of $A$ on $\K$.  Suppose that $V$ is a partial isometry
operating on $\K$ such that for all $x \in I$
\medskip
\itemitem{(i)} $V \pi (x) = \pi (\theta (x)) V$
\medskip
\itemitem{(ii)} $V^* V \pi (x) = \pi (x)$
\medskip
\noindent Then there exists a representation $\pi_\E$ of $\E$ on $\K$
such that
$$\pi_\E (a_n u^n \otimes S^n) = \pi (a_n) V^n \for n \in {\bf Z},
a_n \in D_n.$$
\noindent If, in addition, for all $y \in J$
\medskip
\itemitem{(iii)} $VV^* \pi (y) = \pi (y)$
\medskip
\noindent then there exists a representation $\pi_B$ of $B$ on $\K$
such that
$$\pi_B (a_n u^n) = \pi (a_n) V^n \for n \in \Z, a_n \in D_n.$$}

\proof The second assertion, by far the easier, can be proved
exactly as in the proof of (5.5).  Now assume that $V$ satisfies (i)
and (ii).  Our strategy will be to define a representation $\pi_0$ of
$\E_0$ on $\K$, such that the pair $(\pi_0, V)$ satisfies (i), (ii)
and (iii) with respect to the partial automorphism $(\theta_\E, \E^*_1
\E_1, \E_1 \E^*_1)$ of $\E_0$.  Once that is done, we may apply the
part of the present Lemma which we have already taken care of.
For $n \in \N^*$ define
$$\pi_n \colon  a_n \in D_n \longmapsto V^{n-1} \pi (\theta^{-(n-1)}
(a_n)) V^{{n-1}^*} - V^n \pi(\theta^{-n} (a_n)) V^{n^*} \in B(\K).$$
The reader may now verify that $\pi_n$ is a representation of $D_n$
on $\K$ and that
$$\pi_n (a_n) \pi_m(a_m) = 0$$
for $n \not = m, a_n \in D_n$ and $a_m \in D_m$. In doing so, the
following identities are helpful
$$V^{n^*} V^n \pi (a_{-n}) = \pi (a_{-n}) V^{n^*} V^n = \pi (a_{-n})
\for \quad n \geq 1,\ a_{-n} \in D_{-n }$$
$$\pi (a_n) V^n V^{n^*} V^n = \pi (a_n) V^n \for \quad n \geq 1,\
a_n \in D_n.$$

Next, we may define a representation $\tilde{\pi} $ of $\E_0$ by
$$\tilde{\pi} (a\otimes 1 + \sum_{n=1}^N a_n \otimes e_{nn}) = \pi
(a) + \sum_{n=1}^N \pi_n (a_n).$$
We claim that $(\tilde{\pi}, V)$ satisfy (i), (ii) and (iii) with
respect to $(\theta_\E, \E_1^* \E_1, \E_1 \E^*_1)$. In fact let $x =
a_{-1} \otimes 1 + \sum_{n=1}^N h_n \otimes e_{nn} \in \E_1^* \E_1$
where $a_{-1} \in D_{-1} $ and $h_n \in D_{-1} \cap D_n$.  We have
$$ \tilde{\pi} (\theta_\E (x))V = \tilde{\pi} ((u\otimes S)(a_{-1}
\otimes 1 + \sum_{n=1}^N h_n \otimes e_{nn} ) (u \otimes S)^*) V = $$
$$ \tilde{\pi} ( \theta(a_{-1}) \otimes Q + \sum_{n=1}^N \theta
(h_n) \otimes e_{n+1, n+1} ) V =$$
$$ \tilde{\pi} (\theta (a_{-1}) \otimes 1 - \theta (a_{-1}) \otimes
e_{11} + \sum_{n=1}^N \theta (h_n) \otimes e_{n+1, n+1}) V = $$
$$ (\pi (\theta(a_{-1})) - \pi (\theta (a_{-1})) + V \pi (a_{-1})
V^* + $$
$$ + \sum_{n=1}^N V^n \pi (\theta^{-n+1} (h_n)) V^{n^*} - V^{n+1}
\pi (\theta^{-n} (h_n)) V^{{n+1}^*}) V =$$
 $$ V \pi (a_{-1}) + \biggl( \sum_{n=1}^N V^n V^{{n-1}^*} \pi
(h_n) V^* - V^{n+1} V^{n^*} \pi (h_n) V^* \biggr) V = $$
 $$ V \pi (a_{-1}) + \sum_{n=1}^N V^n V^{{n-1}^*} \pi (h_n) -
V^{n+1} V^{n^*} \pi (h_n) = $$
 $$ V \pi (a_{-1}) + V \sum_{n=1}^N V^{n-1} \pi (\theta^{-(n-1)}
(h_n)) V^{{n-1}^*} - V^n \pi (\theta^{-n} (h_n)) V^{n^*} = $$
 $$ V \tilde{\pi} (a_{-1} \otimes 1 + \sum_{n=1}^N h_n \otimes
e_{nn} ) = V \tilde{\pi} (x).$$
showing (i).  To check (iii) let
$ y = a_1 \otimes Q + \sum^N_{n=2} a_n \otimes e_{nn}$
be in $\E_1 \E^*_1$.  Using that $VV^* V = V$ we have
$$ VV^* \tilde{\pi} (y) = VV^* \tilde{\pi} (a_1 \otimes 1 - a_1
\otimes e_{11} + \sum^N_{n=2} a_n \otimes e_{nn}) = $$
$$ VV^* (V \pi (\theta^{-1} (a_1)) V^* + \sum^N_{n=2} V^{n-1} \pi
(\theta^{-(n-1)} (a_n)) V^{n^*} - V^n \pi (\theta^{-n} (a_n)) V^{n^*})
= \pi (y).$$

The verification of (ii) is left to the reader.  This implies, by
our previous work, that there exists a representation $\pi_\E$ of $\E$
on $\K$ such that
$$\pi_\E (y_n w^n) = \tilde{\pi} (y_n) V^n$$
for $y_n \in \E_n \E^*_n$, where $w$ is the partial isometry
implementing the partial automorphism $\theta_\E$.  But, as it can be
easily verified, $w = u \otimes S$.  Therefore, if $a_1 \in D_1$ we
have that $a_1 \otimes Q \in \E_1 \E_1^*$ and so
$$\pi_\E (a_1 u \otimes S) = \pi_\E ((a_1 \otimes Q) (u\otimes S)) =
\tilde{\pi}(a_1 \otimes Q) V=$$
 $$ \tilde{\pi} (a_1 \otimes 1 - a_1 \otimes e_{11} ) V = V\pi
(\theta^{-1} (a_1)) V^* V = \pi (a_1) V.$$
To conclude the proof one has, for $a_n, b_n \in D_n$ and $n > 1 $,
$$a_n b_n u^n \otimes S^n = (a_n u^{n-1} \otimes S^{n-1} ) (
\theta^{-(n-1)} (b_n) u \otimes S)$$
so, by induction
$$\pi_\E (a_n b_n u^n \otimes S^n) = \pi (a_n) V^{n-1} \pi
(\theta^{-(n-1)} (b_n)) V = \pi (a_n b_n) V^n.$$
\endproof

\statement{6.8. Proposition} {There exists a short exact sequence
$$0 \longrightarrow \Lambda \stackrel{i}{\longrightarrow} \E
\stackrel{\phi}{\longrightarrow} B \longrightarrow 0$$
where $i$ is the inclusion and $\phi$ is given by
$$\phi (a_n u^n \otimes S^n) = a_n u^n \for n \in \Z,\ a_n \in
D_n.$$}

\proof By (6.7) there indeed exists $\phi \colon  \E \rightarrow B$
satisfying $\phi (a_n u^n \otimes S^n) = a_n u^n$.  The identity used
in our proof of (6.4) shows that $\phi$ vanishes on $\Lambda$ so it
defines a homomorphism $\tilde{\phi} \colon  \E/\Lambda \rightarrow B$ which
is obviously surjective.  The action $\gamma$ drops to $\E/\Lambda$,
since $\Lambda$ is an invariant ideal, and then $\tilde{\phi}$ becomes
a covariant homomorphism (with the dual action on $B$).

For covariant homomorphisms under actions of $\S1$, injectivity is
equivalent to injectivity on the fixed point subalgebra by (2.9).  The
latter being true for $\tilde{\phi}$, we conclude that $\tilde{\phi}$
is an isomorphism and the proof is thus completed.  \endproof

\beginsection 7. The Generalized Pimsner-Voiculescu Exact Sequence

In this section we shall obtain a generalization of the
Pimsner-Voiculescu exact sequence \cite{\PV} to our context of crossed
products by partial automorphisms.

An important tool in the sequel will be the fact that the inclusion
of a full corner of a \cstar-algebra induces an isomorphism on
$K$-theory, in the presence of strictly positive elements
(\cite{\Brown}, see also \cite{\PaschkeCircle}, Proposition 2.1). In
fact a much stronger result holds: If $B$ is a full corner in the
\cstar-algebra $A$ and if $A$ has a strictly positive element, then
the inclusion
$\iota \colon  B \rightarrow A$
induces an invertible element in $\KK(B,A)$. This is a trivial
consequence of (\cite{\Brown}, Lemma 2.5). In fact, if $v$ is an
isometry in $M(A\otimes K)$ such that $vv^* = p\otimes 1$, where $p$
is the projection in $M(A\otimes K)$ for which $pAp = B$, then let
$\phi \colon  A\otimes K \rightarrow B\otimes K$
be the isomorphism given by $\phi(x) = v x v^*$. The composition
$\phi (\iota \otimes id_K)$, once composed with the inclusion of
$B\otimes K$ into $M_2(B\otimes K)$, can be easily shown to be
homotopic to the latter map. Likewise, $(\iota \otimes id_K) \phi$
when composed with inclusion of $A\otimes K$ into $M_2(A\otimes K)$ is
also homotopic to the latter map. In a word, $\phi$ provides an
inverse for $\iota$ in $\KK(B,A)$

The model of \KK-theory we shall adopt is the one introduced by
Cuntz in \cite{\CuntzNewLook}.  See also the exposition in Chapter (5)
of \cite{\Jensen} which is where we borrow our notation from.

As we already mentioned, we shall make intensive use of the result
above and hence we need to ensure that several algebras in our
construction, including $\Lambda$, have strictly positive elements.
It is not hard to see that, in the case of $\Lambda$, this implies
that $D_n$ must have a strictly positive element for each $n$. However
we cannot think of a reasonably graceful set of hypothesis which could
provide for this much, other than, of course, assuming our algebras to
be separable (see 1.4.3 and 3.10.5 in \cite{\Pedersen}).

As in the previous section, let $A$ be a \cstar-algebra with a fixed
partial automorphism $\Theta = (\theta, I, J)$ and let $B =
C^*(A,\Theta)$.  In most of our results below, $A$ will, therefore, be
assumed separable.

In the following, all occurrences of $i$ refer to the inclusion
homomorphism that should be clear from the context. Also we shall
denote by $\theta^{-1}_*$ the map induced at the level of $K$ groups
by the composition $$J \buildrel {\theta^{-1}} \over \rightarrow I
\buildrel {i} \over \rightarrow A.$$

The following is our second main result.

\statement{7.1. Theorem} {Let $\Theta = ( \theta, I, J)$ be a
partial automorphism of the separable \cstar-algebra $A$.  Then there
exists an exact sequence of $K$-groups
$$ \matrix {
K_0(J) &\buildrel {i_* - \theta^{-1}_*} \over \longrightarrow
&K_0(A) &\buildrel {i_*} \over \longrightarrow
&K_0(C^*(A,\Theta)) \cr
\cr
\uparrow & & & &\downarrow \cr
\cr
K_1(C^*(A,\Theta)) &\buildrel {i_*} \over \longleftarrow
&K_1(A) &\buildrel {i_* - \theta^{-1}_*} \over \longleftarrow
&K_1(J)
}$$}

\smallskip Our proof will be, roughly speaking, based on \cite{\PV} in
the sense that we shall derive our result from the usual exact
sequence of $K$-Theory for the extension
$$0 \longrightarrow \Lambda \stackrel{i}{\longrightarrow} \E
\stackrel{\phi}{\longrightarrow} C^*(A,\Theta) \longrightarrow 0$$

The crucial part in accomplishing this program is proving that
$K_*(A)$ and $K_* (\E)$ are isomorphic, at which point we can no
longer follow the original method of Pimsner and Voiculescu but,
instead, we must resort to techniques from \KK-theory (see Proposition
5.5 in \cite{\CuntzKK} and Theorem 7.2 in \cite{\AP}) which will, in
fact, lead us to the stronger result that $A$ and $\E$ are
\KK-equivalent to each other.

In the following we let
$$d\colon  A \rightarrow \E$$
and
$$j\colon  J \rightarrow \Lambda $$
be defined by $d(a) = a \otimes 1 $ for $a$ in $A$ and $j(x) = x
\otimes e_{11}$ for $x \in J$.

\statement{7.2. Proposition} {The diagram
$$\matrix{
  K_*(\Lambda) & \stackrel{i_*}{\longrightarrow} & K_*(\E)\cr
  \cr
  j_*\uparrow \ \ \ \ & & \ \ \ \ \uparrow d_*\cr
  \cr
  K_*(J) & \stackrel{i_*- \theta^{-1}_*}{\longrightarrow} & K_*(A)
  }$$
is commutative.  Moreover, if $A$ is separable then $j$ induces a
\KK-equivalence.}

\proof Let $\tilde{J}$ denote $J$ with an added unit and let $y = 1
+ ab$ be an invertible element in $\tilde{J}$ with $a$ and $b$ in $J$.
The $K_1$-class of $y$ being denoted by $[y]_1$, we have
$$d_* \theta^{-1}_* [y]_1 = [1 + u^* abu \otimes 1]_1 = [1+ (u^* a
\otimes S^*) (bu \otimes S)]_1.$$

Recall from Lemma (1.1) in \cite{\PaschkeCircle} that the $K_1$
classes of $1 + rs$ and $1 + sr$ coincide, whenever $1 + rs$ is
invertible.  So we have
$$d_* \theta^{-1}_* [y]_1 = [1 + ba \otimes Q]_1 = [ 1 + (b\otimes
1) (a\otimes Q)]_1 =$$
$$[1 + (a \otimes Q)(b \otimes 1)]_1 = [1 + ab \otimes Q].$$

\noindent Following our diagram in the counterclockwise direction
from $K_1(J)$ we then obtain
$$d_*(i_*-\theta^{-1}_*)[y]_1 = [1 + ab \otimes 1]_1 - [1 + ab
\otimes Q]_1 = [ 1 + ab \otimes e_{11} ]_1 = i_* j_* [y]_1 .$$

Clearly, any invertible element in $\tilde{J}$ is homotopic to an
invertible element of the form $1 + ab$ as above.  Moreover, by
tensoring everything with $M_n ({\bf C})$ we get the above equality
for any invertible $y \in M_n (\tilde J)$.  The commutativity of our
diagram is therefore proved for $K_1$.  The proof for $K_0$ follows by
taking suspensions.

 To conclude, we need to show that $j$ induces a \KK-equivalence.
Note that $j$ is an isomorphism from $J$ onto $J\otimes e_{11} = B_1
B_1^* \otimes e_{11}$.  The latter is obviously a corner of $\Lambda$
and because, for $i, j > 1$
$$(B_i B_1^* \otimes e_{i1}) (B_1 B_1^* \otimes e_{11}) (B_1 B_j^*
\otimes e_{1j}) = $$
$$ B_i B_1^* B_1 B_1^* B_1 B_j^* \otimes e_{ij} = B_i B_j^* \otimes
e_{ij} $$
one sees that $J\otimes e_{11}$ is not contained in any proper ideal
of $\Lambda$.  In two words, $J\otimes e_{11}$ is a full corner of
$\Lambda $.  If we now invoke the separability of $A$, and hence also
of $\Lambda$, the result will follow from our comments in the
introduction to the present section. \endproof

\smallskip It follows that $j_*$ is an isomorphism between the indicated
$K$-groups and so all we need, in order to prove (7.1), is to show
that $d_*$ is an isomorphism as well.  Then we just have to write down
the $K$-theory exact sequence for
$$0 \longrightarrow \Lambda \stackrel{i}{\longrightarrow} \E
\stackrel{\phi}{\longrightarrow} B \longrightarrow 0 .$$
and to replace the two segments $K_* (\Lambda)
\stackrel{i_*}{\longrightarrow} K_* (\E)$ there, by $K_* (J)
\stackrel{i_* - \theta^{-1}_*}{\longrightarrow} K_* (A) $.

\statement{7.3. Definition} {Let $\Lambda_0$ be the subalgebra of
$B(\H\otimes \ell_2 (\N))$ given by
$$\Lambda_0 = \bigoplus_{i,j\geq 0} B_i B_j^* \otimes e_{ij}.$$
(note that, as opposed to (6.3), the indices $i$ and $j$ here, start
at zero). Also let
$$j_0\colon  A \longrightarrow \Lambda_0$$
be given by $j_0 (a) = a \otimes e_{00}$.}

\statement{7.4. Proposition} {If $A$ is separable then $j_0$ is a
\KK-equivalence.}

\proof As in (7.2) it is clear that the image of $j$, namely
$A\otimes e_{00} = B_0 B_0^* \otimes e_{00}$, is a corner of
$\Lambda_0$.  To show that it is also full note that for $i,j \geq 0$
$$(B_i B_0^* \otimes e_{i0} ) (B_0 B_0^* \otimes e_{00} ) (B_0 B^*_j
\otimes e_{0j} ) = B_i B_j^* \otimes e_{ij}.  $$
The result now follows as in (7.2). \endproof

\smallskip Aiming at a proof that $d$ induces an invertible element in
$\KK(A,\E)$, we shall now introduce an element of $\KK (\E, \Lambda_0)$
which will be shown to provide the required inverse, given the
\KK-equivalence between $A$ and $\Lambda_0$.

In order to avoid confusion, we let $T$ be the unilateral shift on
$\ell_2 (\N)$ (as opposed to $S$, which operates on $\ell_2 ({\bf
N}^*))$.  Define
$$\phi, \overline{\phi} \colon \E \longrightarrow B( \H\otimes \ell_2
(\N))$$
by
$$\phi( b_n \otimes S^n) = b_n \otimes T^n$$
and
$$\overline{\phi} (b_n \otimes S^n) = b_n \otimes T^nTT^*.$$
for $n \in \Z$ and $b_n \in B_n$.

\statement{7.5. Proposition} {The images of $\phi$ and
$\overline{\phi}$ are contained in the multiplier algebra
$M(\Lambda_0)$.  Moreover, for every $x\in \E$ we have $\phi (x) -
\overline{\phi}(x) \in \Lambda_0$.}

\proof We have, for $i,j \geq 0$, $n \in \Z$ and $b_n \in B_n$
$$(b_n \otimes T^n) (B_i B_j^* \otimes e_{ij}) = [n+i \geq 0 ]\ b_n
B_i B_j^* \otimes e_{n+i, j}$$
where $[n+i\geq 0]$ indicates the obvious boolean function assuming
the real values 1 and 0.  Also
$$(B_i B_j^* \otimes e_{ij} ) (b_n \otimes T^n) = [j-n \geq 0]\ B_i
B^*_j b_n \otimes e_{i, j-n}$$
showing that $b_n \otimes T^n \in M (\Lambda_0)$.  Similarly the
image of $\overline{\phi}$ is also seen to be contained in
$M(\Lambda_0)$.

Note that for $n \geq 0$, $b_n \in B_n$ and $x = b_n \otimes S^n$
$$\phi (x) - \overline{\phi} (x) = b_n \otimes e_{n,0}$$
which belongs to $\Lambda_0$.  Since the set of $x$'s as above,
generate $\E$, it follows that $\phi(x) - \overline{\phi} (x)$ is in
$\Lambda_0$ for all $x \in \E$. \endproof

\smallskip At this point we start our main \KK-theory computations.  We
refer the reader to chapter (5) in \cite{\Jensen} for the main facts
about the \KK-product as well as for notation.

According to (5.1.1) in \cite{\Jensen}, $q(\phi, \overline{\phi})$
defines a homomorphism from $q \E$ to $\Lambda_0$ and therefore its
homotopy class gives an element, denoted $[q (\phi,
\overline{\phi})]$, in $\KK (\E, \Lambda_0)$.

\statement{7.6. Proposition} {The product $[q (d,0) ] \cdot [q
(\phi, \overline{\phi})]$, in $\KK (A, \Lambda_0)$, equals $[q (j_0,
0)]$.}

\proof Still using the notation from \cite{\Jensen} and thus
denoting by $1_\E$ the unit in the ring $\KK(\E, \E)$, we have
$d^*(1_\E) = [q(d,0)]$. So, by Theorem (5.1.15) in \cite{\Jensen}
$$[q(d,0) ] \cdot [q (\phi, \overline{\phi})] = d^* (1_\E) \cdot [ q
(\phi, \overline{\phi})] = d^* (1_\E \cdot [q (\phi,
\overline{\phi})]) = d^* ([q(\phi, \overline{\phi})]) =$$
$$ [ q (\phi, \overline{\phi}) q(d) ] = [q (\phi d, \overline{\phi}
d )] = [q ( \overline{\phi} d + j_0, \overline{\phi} d)] = [q (j_0,
0)].$$ \endproof

\smallskip Since $[q (j_0, 0)]$ is invertible in $\KK(A, \Lambda_0)$ by
(7.4), if we assume $A$ to be separable, we find that $[q (d,0)]$ is
right invertible and we shall now concentrate our efforts in proving
that it is also left invertible.

\statement{7.7. Definition} {Let $\Omega_0$ be the subset of
$B(\H\otimes \ell_2 (\N^*) \otimes \ell_2 (\N))$ given by
$$\Omega_0 = \bigoplus_{i,j\geq 0} (( B_i \otimes 1) \E (B_j^*
\otimes 1)) \otimes e_{ij}$$}

\smallskip Similarly to what we said with respect to both $\Lambda$ and
$\Lambda_0$, we have that $\Omega_0$ is a subalgebra of operators on
$\H\otimes \ell_2 (\N^*) \otimes \ell_2 (\N)$ while the map
$$k_0 \colon  \E \longrightarrow \Omega_0$$
defined by $k_0 (x) = x \otimes e_{00}$, induces the invertible
element $[q (k_0, 0)]$ in $\KK (\E, \Omega_0)$, when $A$ is separable.

Let $d'$ be the homomorphism
$$d'\colon  \Lambda_0 \longrightarrow \Omega_0$$
defined by
$$d' (b_i b_j^* \otimes e_{ij}) = b_i b_j^* \otimes 1 \otimes e_{ij}
\for i,j \geq 0, \ b_i b_j^* \in B_i B_j^*.$$
Note that the image of $d'$ in fact lies in $\Omega_0$ since, by
(2.7)
$$B_i B_j^* \otimes 1 = B_i B_i^* B_i B_j^* B_j B_j^* \otimes 1 =
(B_i \otimes 1) (B_i^* B_i B_j^* B_j \otimes 1) (B_j^* \otimes 1)
\subseteq$$
$$\subseteq (B_i \otimes 1) (B_0 \otimes 1) (B_j^* \otimes 1)
\subseteq (B_i \otimes 1) \E (B_j^* \otimes 1).  $$

We then have the commutative diagram
$$\matrix{
  \Lambda_0 & \stackrel{d'}{\longrightarrow} & \Omega_0 \cr
  \cr
  j_0 \uparrow \ \ & & \ \ \uparrow k_0 \cr
  \cr
  A & \stackrel{d}{\longrightarrow} & \E
  }$$
in which the vertical arrows are invertible in the corresponding
\KK-groups in the separable case.  It is then obvious that the
existence of a left inverse to $d$ can be verified by exhibiting a
left inverse for $d'$.

\statement{7.8. Proposition} {The product $[q(\phi,
\overline{\phi})] \cdot [q(d',0)]$, in $\KK(\E,\Omega_0)$, equals
$[q (k_0, 0)]$.}

\proof Using (5.1.15) in \cite{\Jensen}, once more, as well as the
fact that $d'_* (1_{\Lambda_{0}}) = [q(d',0)]$, we have
$$[ q (\phi, \overline{\phi})] \cdot [q (d', 0)] = [q(\phi,
\overline{\phi})] \cdot d'_* (1_{\Lambda_{0}}) = d'_* ([q (\phi,
\overline{\phi})] \cdot 1_{\Lambda_{0}}) = d'_* ([q (\phi,
\overline{\phi})])= [ q (\phi, \overline{\phi}) d' ] .$$

Note that, in the last term above, we have written $d'$ where the
definition (see \cite{\Jensen}) of $d'_*$ calls for $id_K \otimes d'$.
We are justified in doing so because $q (\phi, \overline{\phi})$ takes
values in $\Lambda_0$ which should be thought of as a subalgebra of
$K\otimes \Lambda_0$, as is customary in $K$-theory.

Observing that $d'$ extends to a homomorphism from $B(\H\otimes
\ell_2 (\N))$ to $B(\H\otimes \ell_2 (\N^*) \otimes \ell_2 (\N))$ in
an obvious way, one sees that $q (\phi, \overline{\phi}) d' = q (\psi,
\overline{\psi})$ where $\psi$ and $\overline{\psi}$ are the maps from
$\E$ to the multiplier algebra $M(\Omega_0)$ given by
$$\psi (b_n \otimes S^n) = b_n \otimes 1 \otimes T^n$$
and
$$\overline{\psi} (b_n \otimes S^n) = b_n \otimes 1 \otimes T^n
TT^*$$
for $n \in \Z$ and $b_n \in B_n$.  Recalling that $Q = SS^*$ and $P
= 1 - SS^*$, consider the operator $W$ on $\H\otimes \ell_2 (\N^*)
\otimes \ell_2 (\N)$ given by
$$W = 1 \otimes P \otimes e_{00} + 1 \otimes S \otimes e_{01} + 1
\otimes S^* \otimes e_{10} + 1 \otimes 1 \otimes (1 - e_{00} -
e_{11}).$$
It's matrix with respect to the canonic basis of $\ell_2 (\N)$ is
given by
$$W = \pmatrix{
1 \otimes P & 1 \otimes S & 0 & & \cdots \cr
1 \otimes S^* & 0 & 0 & & \cdots\cr
0 & 0 & 1 & & \cdots \cr
& & & 1 \cr
\vdots & \vdots & \vdots & & \ddots
}$$

Clearly, $W$ is a self-adjoint unitary which, therefore can be
connected to the identity operator through the path of unitaries
$$W_t = (1 + W)/2 + e^{\pi i t} (1 - W)/2 \for 0 \leq t \leq 1.$$
If we let $\lambda = \lambda (t) = ( 1 -e^{\pi i t} )/2$ and $\mu =
\mu (t) = (1 + e^{\pi it})/2$ we can write $W_t$ as
$$W_t =
\pmatrix{1 \otimes P + \mu \otimes (1-P) & \lambda \otimes S & 0 & &
\cdots\cr
\lambda \otimes S^* & \mu \otimes 1 & 0 & & \cdots \cr
0 & 0 & 1 & & \cdots \cr
& & & 1 \cr
\vdots & \vdots & \vdots & & \ddots
}$$
Let $V_t = W_t (u\otimes 1 \otimes T)$ or, in matrix form,
$$V_t = \pmatrix{
\lambda u \otimes S & 0 &   & \cdots \cr
\mu u \otimes 1     & 0 &   & \cdots \cr
0     & u \otimes 1     & 0 & \cdots \cr
      & 0 & u\otimes 1      & \ddots \cr
\vdots & \vdots & \ddots    & \ddots \cr
}$$
Observe that, because $W_t$ commutes with $A\otimes 1 \otimes 1$, we
have for all $x \in I$
\medskip
\itemitem{(i)} $V_t (x \otimes 1 \otimes 1) = (\theta (x) \otimes 1
\otimes 1 ) V_t$
\medskip
\itemitem{(ii)} $V_t^* V_t (x \otimes 1 \otimes 1) = x \otimes 1
\otimes 1$.
\medskip
\noindent These are precisely the assumptions in (6.7) so there
exists a representation $\psi_t $ of $\E$ on $\H\otimes \ell_2 (\N^*)
\otimes \ell_2 (\N)$ such that for $n \in \Z$ and $a_n \in D_n$
$$\psi_t (a_n u^n \otimes S^n) = (a_n \otimes 1 \otimes 1) V^n_t$$
Since $V_0 = u \otimes 1 \otimes T$, it is clear that $\psi_0 =
\psi$.  On the other hand it is easy to see that
$$\psi_1 (a_n u^n \otimes S^n) = a_n u^n \otimes S^n \otimes e_{00}
+ a_n u^n \otimes 1 \otimes T^n TT^* \for a_n \in D_n$$
which gives $\psi_1 (x) = k_0 (x) + \overline{\psi} (x) $ for all $x
\in\E$.  Therefore, the pair $ \psi_1, \overline{\psi}$ defines a
homomorphism
$$q (\psi_1, \overline{\psi}) \colon  q \E \longrightarrow \Omega_0$$
which obviously coincide with $q (k_0, 0)$.  To conclude, we need to
show that $\psi_t (x) \in M (\Omega_0)$ and that $ \psi_t (x) -
\overline{\psi} (x) \in \Omega_0$ for all $t$ and for all $x \in \E$.
In fact, once this is done, we obtain a well defined family of
homomorphisms
$$q(\psi_t, \overline{\psi})\colon  q \E \rightarrow \Omega_0,$$
providing a homotopy from $q (\psi, \overline{\psi})$ to $q (k_0,
0)$ and showing that $[q (\psi, \overline{\psi})]$ equals $[q (k_0,
0)]$ in $\KK (\E, \Omega_0)$.

The fact that $\psi_t (x) - \overline{\psi}(x) \in \Omega_0$ is
obvious for $x$ of the form $x = a \otimes 1$ with $a \in A$.  If $x =
a_1 u \otimes S$, for $a_1 \in D_1$, we have
$$\psi_t (x) - \overline{\psi} (x) = \psi_t (x) - \psi (x) + \psi
(x) - \overline{\psi} (x)$$
while
$$\psi_t (x) - \psi (x) = (a_1 \otimes 1 \otimes 1) W_t (u \otimes 1
\otimes T) - (a_1 u \otimes 1 \otimes T) =$$
\medskip
$$ (W_t - I) (a_1 u \otimes 1 \otimes T) =
\pmatrix{\lambda a_1 u \otimes S & 0 & \cdots\cr
-\lambda a_1 u \otimes 1 & 0 & \cdots \cr
0 & 0 & \cdots \cr
\vdots & \vdots & \ddots
}$$
\medskip
which clearly belongs to $\Omega_0$.  This implies two important
facts.  First, since we already know that $\psi (\E) \subseteq M
(\Omega_0)$, we conclude that $\psi_t (x) \in M (\Omega_0)$ for $x \in
(B_0 \otimes 1) \cup (B_1 \otimes S)$, which is a generating set.  So
$\psi_t (x) \in M (\Omega_0)$ for all $x \in \E$.  Secondly, because
$\Omega_0$ is an ideal in $M(\Omega_0)$, our task of proving that
$\psi_t (x) - \overline{\psi} (x) \in\Omega_0$ needs only be verified
on a generating set, which the above computation also yields.
\endproof

\medskip \noindent {\bf Proof of (7.1).} Our last proposition shows
that $[q(d^\prime,0)]$ is left invertible since $[q(k_0 ,0)]$ is
invertible. We conclude that $[q(d,0)]$ is invertible and hence that
$A$ and $\E$ are \KK-equivalent to each other. In particular the map
$$d_* \colon  K_*(A) \rightarrow K_*(\E)$$
is an isomorphism.  Our result now follows from combining the
$K$-theory exact sequence for the Toeplitz extension
$$0 \longrightarrow \Lambda \stackrel{i}{\longrightarrow} \E
\stackrel{\phi}{\longrightarrow} B \longrightarrow 0$$
with the conclusion of (7.2).

\bigskip
\bigskip
\centerline{\scX References}
\frenchspacing

\medskip
\stdbib{\AP}
{J. Anderson and W. Paschke}
{The $K$-Theory of the reduced \cstar-algebra of an HNN-group}
{J. Oper. Theory}
{16}
{(1986), 165-187}

\stdbib{\Brown}
{L. G. Brown}
{Stable Isomorphism of Hereditary Subalgebras of \cstar-Algebras}
{Pacific J. Math.}
{71}
{(1977), 335-348}

\stdbib{\BGR}
{L. G. Brown, P. Green, M. A. Rieffel}
{Stable Isomorphism and Strong Morita Equivalence of \cstar-algebras}
{Pacific J. Math.}
{71}
{(1977), 349-363}

\stdbib{\CuntzOnI}
{J. Cuntz}
{$K$-Theory for certain \cstar-Algebras}
{Ann. Math.}
{113}
{(1981), 181--197}

\stdbib{\CuntzOnII}
{\SameAuthor}
{$K$-Theory for certain \cstar-Algebras, II}
{J. Oper. Theory}
{5}
{(1981), 101--108}

\bib{\CuntzKK}
{\SameAuthor}
{$K$-Theory and \cstar-Algebras}
{Lecture Notes in Math. 1046, Springer--Verlag, 1984, pp. 55--79}

\stdbib{\CuntzNewLook}
{\SameAuthor}
{A new look at \KK-Theory}
{$K$-Theory}
{1}
{(1987), 31--51}

\stdbib{\FellI}
{J. M. G. Fell}
{An extension of Mackey's method to Banach *-Algebraic Bundles}
{Memoirs Amer. Math. Soc.}
{90}
{(1969)}

\bib{\FellII}
{\SameAuthor}
{Induced Representations and Banach $*$-Algebraic Bundles}
{Lecture Notes in Math. 582, Springer-Verlag, 1977}

\bib{\Jensen}
{K. Jensen, K. Thomsen}
{Elements of \KK-Theory}
{Birkh\"auser, 1991}

\stdbib{\KishimotoTakai}
{A. Kishimoto, H. Takai}
{Some Remarks on \cstar-Dynamical Systems with a Compact Abelian Group}
{Publ. RIMS, Kyoto Univ.}
{14}
{(1978), 383-397}

\stdbib{\Kumjian}
{A. Kumjian}
{On localizations and simple \cstar-Algebras}
{Pacific J. Math.}
{112}
{(1984), 141-192}

\stdbib{\PaschkeEndom}
{W. L. Paschke}
{The Crossed Product of a \cstar-algebra by an Endomorphism}
{Proc. Amer. Math. Soc.}
{80}
{(1980), 113-118}

\stdbib{\PaschkeCircle}
{\SameAuthor}
{$K$-Theory for actions of the circle group on \cstar-algebras}
{J. Oper. Theory}
{6}
{(1981), 125-133}

\bib{\Pedersen}
{G. K. Pedersen}
{\cstar-Algebras and their Automorphism Groups}
{Academic Press, 1979}

\bib{\Phillips}
{C. Phillips}
{Equivariant $K$-Theory and Freeness of Group Actions on \cstar-algebras}
{Lecture Notes in Math. 1274, Springer-Verlag, 1987}

\stdbib{\PV}
{M. Pimsner and D. Voiculescu}
{Exact sequences for $K$-groups and Ext-groups of certain cross-products
\cstar-algebras}
{J. Oper. Theory}
{4}
{(1980), 93-118}

\stdbib{\RieffelInduced}
{M. A. Rieffel}
{Induced Representations of \cstar-algebras}
{Adv. Math.}
{13}
{(1974), 176-257}

\bib{\RieffelMorita}
{\SameAuthor}
{Morita Equivalence for Operator Algebras}
{Proc. Symp. Pure Math. 38, Amer. Math. Soc., (1982), 285-298}

\bib{\RieffelProper}
{\SameAuthor}
{Proper actions of Groups on \cstar-Algebras}
{Mappings of Operator Algebras, Proc. Japan--US joint seminar, Birkh\"auser
(1990), 141-182}

\bigskip
\rightline{October 1992}

\end